%% file: progressbasisopt_final.tex
\documentclass[times]{nmeauth}
\pdfoutput=1

\usepackage{bm}
\usepackage{amsfonts}
\usepackage{graphicx}
\usepackage{amsthm}
\usepackage{amssymb}
\usepackage{paralist}
\usepackage[section]{placeins}
\usepackage{algorithm}
\usepackage{algorithmic}
\usepackage{color}
\usepackage{moreverb}
\usepackage[colorlinks,bookmarksopen,bookmarksnumbered,citecolor=red,urlcolor=red]{hyperref}
\usepackage{array}
\usepackage{booktabs}
\usepackage{footnote}
\usepackage{caption}
\usepackage{subcaption}

\input{mycommands}

\begin{document}

\runningheads{M.J.~ZAHR and C.~FARHAT}{PDE-CONSTRAINED OPTIMIZATION USING A PROGRESSIVELY CONSTRUCTED ROM}

\title{Progressive construction of a parametric reduced-order model for PDE-constrained optimization}

\author{Matthew J. Zahr \affil{1}\corrauth and Charbel Farhat \affil{1,}\affil{2,}\affil{3}}

\address{\centering \affilnum{1} Institute for Computational and Mathematical Engineering \\
            \affilnum{2} Department of Aeronautics and Astronautics \\
           \affilnum{3} Department of Mechanical Engineering \\
          Stanford University, Mail Code 4035, Stanford, CA 94305, U.S.A.}

\corraddr{Department of Aeronautics and Astronautics, Stanford University, Mail Code 4035, Stanford, CA 94305, U.S.A.}

\begin{abstract}

An adaptive approach to using reduced-order models as surrogates in PDE-constrained optimization is introduced that breaks the traditional offline-online framework of model order reduction. A sequence of optimization problems constrained by a given Reduced-Order Model (ROM) is defined with the goal of converging to the solution of a given PDE-constrained optimization problem.  For each reduced optimization problem, the constraining ROM is trained from sampling the High-Dimensional Model (HDM) at the solution of some of the previous problems in the sequence.  The reduced optimization problems are equipped with a nonlinear trust-region based on a residual error indicator to keep the optimization trajectory in a region of the parameter space where the ROM is accurate.  A technique for incorporating sensitivities into a Reduced-Order Basis (ROB) is also presented, along with a methodology for computing sensitivities of the reduced-order model that minimizes the distance to the corresponding HDM sensitivity, in a suitable norm.  The proposed reduced optimization framework is applied to subsonic aerodynamic shape optimization and shown to reduce the number of queries to the HDM by a factor of 4-5, compared to the optimization problem solved using only the HDM, with errors in the optimal solution far less than $0.1\%$.

\end{abstract}


\maketitle

\section{INTRODUCTION}
Optimization problems constrained by nonlinear Partial Differential Equations (PDE) arise in many engineering applications, including inverse modeling, optimal control, and design.  In the context of design, PDE-constrained optimization provides a valuable tool for optimizing specific features of an existing design.  In practice, the solution of the optimization problem will require many solutions of the underlying PDE, which will be expensive when the underlying PDE model is large.  Large-scale PDE models commonly arise in industry-scale applications where high-fidelity simulations are required due to the complexity of the underlying physics.  In particular, large-scale, nonlinear PDE models arise in Computational Mechanics (CM) simulations, such as viscous or high-speed Computational Fluid Dynamics (CFD), large-deformation Computational Structural Dynamics (CSD) of complex bodies, and high-frequency acoustics, to name a few.

In practice, PDE-constrained optimization problems are solved using a Nested Analysis and Design (NAND)  \cite{biegler2003large} approach, whereby the PDE state variable is considered an implicit function of the parameters of the PDE.  The implication of this approach is that the PDE must be solved at every parameter value encountered in the design cycle.  It is well-known that this approach requires a large number of PDE solutions, and sensitivities in a gradient-based approach, which is prohibitively expensive for industry-scale PDE models.  Accordingly, the main objective of this paper is to present a method for finding an approximate solution to the PDE-constrained optimization problem at a reduced computational cost. This method relies on using a sequence of low-dimensional reduced-order models as surrogates for the underlying high-dimensional model.

Model Order Reduction (MOR) has become an attractive method for reducing the dimension of a PDE by constraining the solution to lie in a low-dimensional affine subspace.  MOR has found the most success in real-time applications, such as optimal control \cite{ito1998reduced,ravindran2000reduced,afanasiev2001adaptive,fahl2001trust,bergmann2005optimal,bergmann2008optimal,kunisch2008proper} and non-destructive evaluation \cite{grepl2007certified}, and many-query applications \cite{amsallem2010towards,fra3}, such as design optimization \cite{legresley2000airfoil,rozza2010model,lassila2010parametric,manzoni2012shape} and uncertainty quantification \cite{buithanh2008parametric}.  Model reduction is often embedded in an offline-online framework, whereby a low-dimensional Reduced-Order Basis (ROB)
is constructed in the \emph{non-time-critical} offline phase.  The span of the ROB defines the subspace in which perturbations of the PDE solution from a chosen offset will be constrained to lie.  Subsequently, the ROB is used to project the governing equations, in either a Galerkin or Petrov-Galerkin sense, onto a low-dimensional test subspace.  The reduced equations and solution space collectively define the 
reduced-order model, which can be used in the \emph{time-critical} online phase.

In this work, projection-based ROMs are considered with bases built using a variant of Proper Orthogonal Decomposition (POD) and the method of snapshots \cite{sirovich1987turbulence}.  In this approach, information is collected from simulations of the HDM at various parameter configurations and compressed to generate the ROB.  The data collected from the HDM simulation are usually state vectors (time history or steady state), but could be any quantity that provides information about the subspace on which the solution lies, such as state sensitivities, residuals, or Krylov vectors \cite{wash2013comm}.  In the remainder, a parameter value at which a HDM is simulated to generate snapshots will be referred to as a \emph{HDM sample}.

A well-known bottleneck in the online phase of projection-based ROMs of nonlinear systems is evaluation and projection of high-dimensional nonlinear quantities that cannot be precomputed in the offline phase.  Only in the case of polynomial nonlinearities can all high-dimensional operations be confined to the offline phase.  In the general case, despite the reduced dimensionality of the model, in the sense of a reduced set of equations and unknowns, the computation time is not necessarily reduced, just redistributed.  Hyperreduction, a term introduced in \cite{ryckelynck2005priori}, has been developed as an additional level of approximation built on top of the ROM to realize computational speedups in addition to dimensionality reduction.  Many different flavors of hyperreduction exist \cite{barrault2004empirical,chaturantabut2010nonlinear,carlberg2011gnat,carlberg2013gnat,ECSW1}, some of which are variants of Gappy POD \cite{everson1995karhunen} and have been shown to demonstrate promising (CPU time) speedup factors in excess of 400 in computing the unsteady solution of a large-scale CFD problem \cite{carlberg2013gnat}.  The method proposed in this work will only be applied at the level of the ROM; extensions to hyperreduction will be mentioned throughout, with complete formulation deferred to future work.  

The standard offline-online \cite{legresley2000airfoil,oliveira2007reduced,bui2008model,rozza2010model,lassila2010parametric,manzoni2012shape,urban2012greedy,hinze2013model} decomposition of computational effort for reduced-order modeling is most appropriate in the real-time application context, where there is a natural offline-online breakdown and minimizing online computing time is paramount.  In the context of many-query analyses, particularly design optimization, the offline-online decomposition is artificial in that the only interest is total wall-time to obtain an accurate solution.  This implies the existence is a ``break-even point'' \cite{huynh2007reduced,oliveira2007reduced,carlberg2008compact,rozza2008reduced,carlberg2011low}, an upper limit on the offline time, beyond which, it cannot be amortized in the online phase.  Furthermore, it is well-known that ROMs lack robustness with respect to parametric variations \cite{epureanu03_ParamROMinTurbo,lieu2006reduced,lieu2007adaptation,hay2009local,amsallem2008interpolation,amsallem2009method}.  
The assumption that the PDE solution can be represented by a single, static, low-dimensional subspace will break-down, particularly for nonlinear problems where the solution can change dramatically with parameter perturbations, such as high-speed or turbulent (chaotic) fluid flows.  As the trajectory of the optimization iterates is unknown a-priori, a ROB must be trained such that it maintains accuracy throughout the parameter space if the offline-online decomposition is used.  The implication is that a large ROB built from extensive training with the HDM will be required.    Despite the extensive training, only the portion of the HDM samples near the optimization trajectory will be relevant.

Reduced-order models have been used as surrogates in PDE-constrained optimization in two main contexts: offline-online framework \cite{legresley2000airfoil,oliveira2007reduced,bui2008model,rozza2010model,manzoni2012shape,urban2012greedy,hinze2013model,dihlmann2013certified}, where efficient HDM sampling \cite{bui2008model,urban2012greedy,hinze2013model} and ROB construction is paramount and an adaptive framework that merges the offline and online phases \cite{ravindran2000reduced,afanasiev2001adaptive,fahl2001trust,bergmann2005optimal,bergmann2008optimal,kunisch2008proper,yue2013accelerating,zahr2013construction}.  The adaptive approaches to ROM-constrained optimization involve a sequence of reduced optimization problems, where the solution to a given reduced optimization problem is used to enrich the training, resulting in an improved basis for subsequent reduced optimization problems.  Variations include trust region methods that treat the ROM objective as the model problem instead of the traditional quadratic model problem -- Trust-Region Proper Orthogonal Decomposition (TRPOD) \cite{fahl2001trust}, using the gradient of the HDM to take a single optimization step before enriching the ROB -- Optimality System Proper Orthogonal Decomposition (OS-POD) \cite{kunisch2008proper}, and equipping the reduced optimization problem with a constraint (or penalty) on the error in the ROM, in the context of linear systems \cite{yue2013accelerating}.  Additionally, in \cite{zahr2013construction}, an adaptive approach was introduced that defined a sequence of reduced optimization problems that balanced minimizing the objective function with targeting locations in the parameter space where the ROM lacks accuracy, i.e. a desirable sample point.

The proposed approach attempts to combine the desirable features of TRPOD, OS-POD, and the error-aware reduced optimization approach for general nonlinear problems, cognizant of the reality that hyperreduction will eventually be required.  The approach defines a sequence of reduced optimization problems, where the ROM at a given point in the sequence is trained from snapshots of the HDM sampled at the solution of previous problems in the sequence. This is a common feature of aforementioned adaptive methods.  Each reduced optimization problem is equipped with a constraint that sets an upper bound on the norm of the HDM residual evaluated at the solution of the ROM.  The norm of the HDM residual is used as an error indicator for the ROM, as efficient and computable error bounds are not available for POD-based model order reduction in the general nonlinear parametric setting.  The HDM residual norm and upper bound collectively define a nonlinear trust region.  The upper bound on the nonlinear trust-region radius will be adaptively refined using an algorithm motivated from standard trust-region methods \cite{fahl2001trust,nocedal2006numerical}.  Therefore, the proposed method incorporates generalizations of the error-awareness of \cite{yue2013accelerating} and the trust region ideology of \cite{fahl2001trust}.
Construction of a ROB and model order reduction framework such that the ROB accurately represents state vectors and sensitivities without degrading the approximation of either is crucial in the context of ROM-constrained optimization \cite{fahl2001trust}.  This issue is addressed in this work by collecting state and sensitivity snapshots separately, applying POD to each snapshot set, and concatenating the results to obtain the ROB \cite{wash2013comm}.  Degradation of the ROM solution and their sensitivities are avoided by applying minimum-residual methods for solving for each.  In particular, the Least-Squares Petrov-Galerkin projection from \cite{barrault2004empirical,bui2008model,carlberg2011gnat} will be used to define the reduced \emph{state} equations, which has a minimum-residual property thereby ensuring that, under some assumptions, the ROM approximation to the PDE state variables can only improve as vectors are added to the ROB.   In the present situation, this implies that the presence of sensitivity snapshots in the ROB cannot adversely affect the accuracy of the state approximation.  Additionally, a novel minimum-error reduced sensitivity framework is introduced that ensures the reduced sensitivity approximation to the PDE sensitivities can only improve as vectors are added to the ROB, i.e. state snapshots in the ROB will not adversely affect the sensitivity approximation.  The minimum-residual formulations for the reduced state and sensitivity equations, along with careful selection of the offset and basis vectors used to define the low-dimensional affine subspace of the ROM, will be shown to generate ROMs whose solution and sensitivities exactly match those of the HDM at training points, provided mild conditions are met. Finally, fast basis updating techniques \cite{wash2012locbaseup} will be discussed for adding snapshot contributions to a ROB without requiring an SVD computation of the entire snapshot matrix appended with new data.

In this paper, \emph{steady-state}, nonlinear PDE constraints are considered using the NAND framework in the discretize-then-differentiate setting for PDE-constrained optimization.  A variant of POD, presented in Section~\ref{sec:offline}, is used to generate an optimization-oriented ROB from state and sensitivity snapshots collected from sampling the HDM at various parameter values.  Although the proposed method for ROM-constrained optimization is derived in this framework, it can be generalized to other PDE-constrained optimization frameworks and other methods of ROB construction.

The remainder of the paper is organized as follows.  Section~\ref{sec:pdeopt} presents a brief overview of PDE-constrained optimization.  Section~\ref{sec:rom} introduces ROMs and a technique for including sensitivity snapshots in the construction of the ROB that will have implications in the context of ROM-constrained optimization.  Section~\ref{sec:progromopt} presents the proposed approach for solving a PDE-constrained optimization problem with a progressively-constructed ROM, including algorithmic details that guarantee desirable properties for each reduced optimization problem.  Section~\ref{sec:app} demonstrates the capabilities of the progressively-constructed ROM approach in solving relevant CFD problems.  Section~\ref{sec:conc} offers conclusions and avenues for further research.

\section{PDE-CONSTRAINED OPTIMIZATION} \label{sec:pdeopt}
In this section, the framework chosen for PDE-constrained optimization is described.  It is emphasized that these choices are made solely due to convenience of integration with currently available tools.  The proposed approach to ROM-constrained optimization in Section~\ref{sec:progromopt} is general in that it is applicable to any PDE-constrained optimization framework.

Recall the focus of this work is \emph{steady-state} partial differential equations.  Define the PDE-constrained optimization problem as

\begin{equation} \label{opt:continuous}
\begin{aligned}
& \underset{\wbold \in \boldsymbol{\Wcal},~\mubold \in \boldsymbol{\Dcal}}{\text{minimize}}
& & \hat\Jcal\left(\wbold, \mubold\right) \\
& \text{subject to}
& & \hat\cbold(\wbold, \mubold) \leq 0 \\
& & &  \boldsymbol{\Lcal}(\wbold, \mubold) = 0,
\end{aligned}
\end{equation} 
where $\boldsymbol{\Wcal}$ is the state space of the corresponding PDE, $\boldsymbol{\Dcal} \subset \Rbb^{n_p}$ is the parameter space, $\func{\hat\Jcal}{\boldsymbol{\Wcal} \times \boldsymbol{\Dcal}}{\Rbb}$ is the objective function, $\func{\hat\cbold}{\boldsymbol\Wcal \times \boldsymbol\Dcal}{\Rbb^{n_c}}$ are $n_c$ additional constraints on the system, and $\boldsymbol{\Lcal}$ is the nonlinear, vector-valued differential operator defining the PDE of interest.  In the following, the discretization of (\ref{opt:continuous}) for integration with an optimization solver is discussed.  The treatment of the PDE constraint and derivation of sensitivities are also discussed.

\subsection{Discretize-then-differentiate vs. differentiate-then-discretize}
The two approaches to deriving a discretized optimization problem of (\ref{opt:continuous}) are known as discretize-then-differentiate and differentiate-then-discretize \cite{gunzburger2003perspectives,pinnau2008optimization}.  In the discretize-then-differentiate approach, the PDE is discretized prior to the derivation of the sensitivity equations, while the differentiate-then-discretize approach derives the sensitivity equations at the PDE-level and discretizes the result.  The operations of discretization and differentiation do not commute and the distinction between these approaches is significant.

The discretize-then-differentiate approach is used throughout this work, however; formulation in the differentiate-then-discretize can be developed similarly.  At this point, the PDE constraint in (\ref{opt:continuous}) is discretized, resulting in the following optimization problem
\begin{equation} \label{opt:hdm}
\begin{aligned}
& \underset{\wbold \in \Rbb^N,~\mubold \in \boldsymbol{\Dcal}}{\text{minimize}}
& & \Jcal\left(\wbold, \mubold\right) \\
& \text{subject to}
& & \cbold(\wbold, \mubold) \leq 0 \\
& & &  \boldsymbol{\Rbold}(\wbold, \mubold) = 0,
\end{aligned}
\end{equation} 
where $\wbold \in \Rbb^N$ is the discretized PDE state vector, $\func{\Jcal}{\Rbb^N \times \boldsymbol{\Dcal}}{\Rbb}$ is the discretized objective function, $\func{\cbold}{\Rbb^N \times \boldsymbol{\Dcal}}{\Rbb^{n_c}}$ are the $n_c$ additional constraints (discretized), and $\Rbold(\wbold,\mubold) = 0$ is the discretized PDE.  In the next section, treatment of the discretized PDE constraint and sensitivity derivation is discussed.

\subsection{Nested analysis and design (NAND)} \label{sec:nand}
The NAND \cite{biegler2003large} approach to PDE-constrained optimization uses the PDE to implicitly define the state vector as a function of the parameters.  In this setting, the only optimization variables are the PDE parameters; however, any parameter value change requires solving the PDE to obtain the corresponding state.  This can be viewed as enforcing satisfaction of the discrete PDE constraint at every iteration of the optimization, whereas most algorithms allow infeasibilities during intermediate iterations.  In general, the NAND approach enables the use of a ``black-box'' PDE solver, although many PDE solutions may be required during the optimization process.  Alternatives to the NAND approach include full-space and reduced-space Simultaneous Analysis and Design (SAND) \cite{biegler2003large,akccelik2006parallel} methods.  In this work, the proposed method is formulated in the NAND framework; however, it is equally
applicable to SAND frameworks.

In the NAND approach, $\wbold = \wbold(\mubold)$ implicitly through the discretized PDE, $\Rbold(\wbold,\mubold) = 0$.  Additionally, satisfaction of the PDE constraint for all $\mubold \in \boldsymbol{\Dcal}$ in the optimization procedure implies that the PDE states lie on the manifold defined by
\begin{equation}\label{eqn:state-man}
  \Rbold(\wbold(\mubold),\mubold) = 0.
\end{equation}
On this manifold, $$\oder{\Rbold}{\mubold} = \pder{\Rbold}{\mubold} + \pder{\Rbold}{\wbold}\pder{\wbold}{\mubold} = 0$$
holds and is used to compute the sensitivities of the state variables with respect to the PDE parameters
\begin{equation} \label{eqn:hdm-sens}
  \pder{\wbold}{\mubold} = -\left[\pder{\Rbold}{\wbold}\right]^{-1}\pder{\Rbold}{\mubold},
\end{equation}
which is subsequently used to compute the gradients of the objective function and constraints.  

Consider some functional $\Fcal \left(\wbold(\mubold), \mubold\right)$, which may represent the objective function or a constraint.  The gradient of $\Fcal$ with respect to $\mubold$ takes the form
\begin{align}
  \oder{\Fcal}{\mubold} = \pder{\Fcal}{\mubold} +  \pder{\Fcal}{\wbold}\pder{\wbold}{\mubold} &= \pder{\Fcal}{\mubold} -\pder{\Fcal}{\wbold}\left(\left[\pder{\Rbold}{\wbold}\right]^{-1}\pder{\Rbold}{\mubold}\right)\label{eqn:sens-direct}\\
  &= \pder{\Fcal}{\mubold} -\left[\pder{\Rbold}{\wbold}^{-T}\pder{\Fcal}{\wbold}^T\right]^T\pder{\Rbold}{\mubold} \label{eqn:sens-adj}.
\end{align}
Equations (\ref{eqn:sens-direct}) and (\ref{eqn:sens-adj}) suggest two ways of computing the gradient of the functional, known as the direct method and adjoint method.  In the direct approach in (\ref{eqn:sens-direct}), the state sensitivity is computed directly, requiring the solution of $n_p$ linear systems with the PDE Jacobian matrix.  This may be very expensive, but is only required once, regardless of the number of functionals that need to be differentiated, allowing for an arbitrary number of objectives and constraints at the cost of $n_p$ linear systems solves.  The adjoint approach in (\ref{eqn:sens-adj}) requires a single linear system solution with the transpose of the PDE Jacobian, regardless of the number of parameters involved.  However, such a solve is required for each functional to be differentiated.  Also notice that $\displaystyle{\pder{\wbold}{\mubold}}$ is never directly computed in the adjoint method.  Therefore, the direct method is preferred over the adjoint method when there are a large number of functionals (objective and constraint functions) compared to the number of parameters and vice versa.  In this work, sensitivities of the state vector, $\displaystyle{\pder{\wbold}{\mubold}}$, are to be used as snapshots in the training of a ROB, see Section~\ref{sec:offline}, thereby motivating the use of the direct method.

\section{MODEL REDUCTION} \label{sec:rom}

\subsection{High-dimensional model}

The discretization of the governing PDE will serve as a set of nonlinear equality constraints in the optimization problem.  The discretized, steady, nonlinear PDE takes the form
  \begin{equation} \label{eqn:hdm-ss}
    \Rbold(\wbold, \mubold) = 0
  \end{equation}
  where $\wbold \in \Rbb^N$ is the discretized state vector, $\mubold \in \boldsymbol{\Dcal} \subset \Rbb^{n_p}$ represents the parameters of the PDE, i.e. shape, control, material, and $N$ assumed to be large.  Depending on the application, it is common practice to use homotopy to solve the above nonlinear function as opposed to applying a nonlinear solver to (\ref{eqn:hdm-ss}) directly.  Since a good initial guess is usually not known for (\ref{eqn:hdm-ss}), the goal of homotopy is to generate a sequence of ``easier'' subproblems, each with an obvious, high-quality initial guess.  Common applications of homotopy are: structural mechanics where load stepping is used and CFD where pseudo-time-stepping is the standard \cite{kelley1998convergence,kelley2007projected}.

\subsection{Reduced-order model} \label{sec:rom-theory}
The goal of model order reduction is to reduce the HDM (\ref{eqn:hdm-ss}), in the sense of reducing the number of degrees of freedom and nonlinear equations, with the expectation of computational savings in terms of both storage and CPU time.  The central assumption in MOR is the state vector lies (approximately) in a low-dimensional affine subspace
\begin{equation}\label{eqn:rom-assump}
  \wbold \approx \wbold_r = \bar\wbold + \Phibold \ybold,
\end{equation}
where $\Phibold \in \Rbb^{N \times k_y}$ is the ROB of dimension $k_y \ll N$.  The offset vector $\bar\wbold$ can be used to make the reduced approximation \emph{exact} at a particular point in state space, namely $\bar\wbold$, usually the initial condition.  This implies that the ROB can be constructed to represent \emph{perturbations} about $\bar\wbold$, instead of the state vector itself \cite{wash2012locbaseup}.  The offset vector $\bar\wbold$ is assumed constant or piecewise constant with respect to $\mubold$ and is revisited in Sections~\ref{sec:offline} and~\ref{sec:progromopt}.  Implicit in assumption (\ref{eqn:rom-assump}) is the sensitivities of the state vector can be represented in \emph{the same} low-dimensional subspace, modulo the affine offset
\begin{equation} \label{eqn:rom-sens-assump}
  \pder{\wbold}{\mubold} \approx \pder{\wbold_r}{\mubold} = \Phibold\pder{\ybold}{\mubold},
\end{equation}
which is obtained by differentiating (\ref{eqn:rom-assump}) with respect to the parameters, $\mubold$.

Substituting the assumption in (\ref{eqn:rom-assump}) into the HDM in (\ref{eqn:hdm-ss}) and constraining the nonlinear equations to be orthogonal to the subspace spanned by a left subspace $\Psibold \in \Rbb^{N \times k_y}$, the ROM is
\begin{equation} \label{eqn:rom}
  \Psibold^T\Rbold(\bar{\wbold} + \Phibold \ybold; \mubold) = 0.
\end{equation}

It was shown in \cite{carlberg2011gnat} that the left basis can be constructed such that the nonlinear search directions of (\ref{eqn:rom}) minimize the distance to the Newton direction of the HDM, in some norm.  For $\Psibold = \Phibold$, (\ref{eqn:rom}) corresponds to a Galerkin projection of the HDM (\ref{eqn:hdm-ss}) and leads to ``optimal'' search directions in the $\displaystyle{\pder{\Rbold}{\wbold}}$ norm for problems with 
Symmetric-Positive Definite (SPD) Jacobians.  For $\displaystyle{\Psibold(\wbold, \mubold) = \pder{\Rbold}{\wbold}(\wbold,\mubold)\Phibold}$, (\ref{eqn:rom}) is a Least-Squares Petrov-Galerkin (LSPG) \cite{barrault2004empirical,bui2008model,carlberg2011gnat} projection, which is ``optimal'' in the $\displaystyle{\pder{\Rbold}{\wbold}^T\pder{\Rbold}{\wbold}}$ norm for problems with Jacobians that are not necessarily SPD.

It can be shown that the LSPG form of (\ref{eqn:rom}) is precisely the first-order optimality condition of 
\begin{equation} \label{eqn:min-res}
	\begin{aligned}
	& \underset{\ybold \in \Rbb^{k_y}}{\text{minimize}}
	& & \frac{1}{2} \norm{\Rbold\left(\bar\wbold + \Phibold\ybold, \mubold\right)}_2^2.
	\end{aligned}
\end{equation}
  A ROM that is equivalent to a minimization problem of the form (\ref{eqn:min-res}) will be said to possess the minimum-residual property.  The minimum-residual property is desirable as it guarantees that additional information in a ROB cannot degrade the solution approximation, neglecting issues related to convergence and assuming a monotonic relationship between residual norm and error.
In this work, the general case of the LSPG ROM is used.  The proposed method is equally applicable to Galerkin ROMs and the derivation of the reduced sensitivities simplifies as the left basis does not depend on the state vector, as will be shown in Section~\ref{sec:rom-sens}.
 
 As previously mentioned, the projection-based ROM in (\ref{eqn:rom}) is not sufficient to generate non-trivial speedups due to the large-dimensional operations that must be done at every nonlinear iteration to evaluate and project the nonlinear terms.  A slew of \emph{hyperreduction} techniques have been developed, most of which aim to approximate the nonlinear terms in a low-dimensional subspace \cite{chaturantabut2010nonlinear,carlberg2011gnat}.  An important component of these hyperreduction techniques is the concept of \emph{masked} or gappy nonlinear terms, defined as a subset of the entries of a high-dimensional nonlinear quantity.  The masked nonlinear terms are used to determine the remaining entries via interpolation or least-squares minimization.  Additionally, the sample mesh is defined as the subset of the entire mesh required to compute the masked nonlinear terms.  This approach has been shown to generate promising speedups as the nonlinear terms only need to be evaluated on a small subset of the original mesh and the large-dimensional projections can be precomputed.
 
 In this work, the proposed method is only derived for the projection-based ROM (\ref{eqn:rom}); however, it is equally applicable in the context of hyperreduction with several complications that arise when only masked nonlinear terms are available.  In the remainder of this paper, extensions to hyperreduction are mentioned, although a full exposition using hyperreduced models as a surrogate for PDE-constrained optimization is postponed to future work.

\subsection{Offline basis construction: proper orthogonal decomposition and the method of snapshots} \label{sec:offline}
The method of snapshots in conjunction with Proper Orthogonal Decomposition (POD) \cite{sirovich1987turbulence} is among the most popular and successful techniques for collecting and compressing information to yield a ROB in the context of model order reduction.  The method of snapshots consists of sampling the HDM at various parameter configurations and collecting snapshots, usually state vectors, from the simulations.  POD, also known as Principal Component Analysis (PCA) or the Karhunen-Loeve procedure is used to compress the snapshots, exploiting the optimal approximation properties of the Singular Value Decomposition (SVD), to yield a ROB.  It is worth mentioning that a complete singular value decomposition is among the most expensive matrix factorizations implying that basis construction via POD can constitute a significant portion of the traditional offline phase, particularly for a large dimensional state space and many snapshots.

\begin{algorithm}[htbp]
\caption{Proper Orthogonal Decomposition}
\begin{algorithmic}[1]\label{alg:pod}
\REQUIRE Snapshot matrix, $\Xbold \in \Rbb^{N \times k_s}$ and ROB size, $k_y$
\ENSURE  ROB, $\Phibold$
  \STATE Compute the thin SVD of $\Xbold$: $\Xbold = \Ubold\Sigmabold\Vbold^T$, where $\Ubold = \begin{bmatrix} \ubold_1 & \ubold_2 & \cdots & \ubold_{k_s} \end{bmatrix}$
  \STATE $\Phibold = \begin{bmatrix} \ubold_1 & \ubold_2 & \cdots & \ubold_{k_y} \end{bmatrix}$
\end{algorithmic}
\end{algorithm}

From the low-dimensionality assumptions in (\ref{eqn:rom-assump}) and (\ref{eqn:rom-sens-assump}), it is clear that the subspace spanned by the ROB needs to (approximately) contain $\wbold(\mubold) - \bar\wbold$ and $\displaystyle{\pder{\wbold}{\mubold}(\mubold)}$ if the ROM is to accurately represent the HDM solution and sensitivity at $\mubold$.  This suggests using a snapshot matrix of the form 
\begin{equation} \label{eqn:snap-cat}
  \Xbold = \left[\wbold(\mubold) - \bar\wbold, \pder{\wbold}{\mubold}(\mubold)\right]
\end{equation}
is desirable in the context of optimization, as recognized in \cite{arian2000trust,carlberg2008compact,carlberg2011low}.  As POD is sensitive to the relative scale of the columns of the data matrix, it cannot be directly applied to a snapshot matrix containing states and sensitivities due to the inherent scale difference.  A heuristic weighting procedure was developed in \cite{carlberg2008compact,carlberg2011low} to address the scale differences between the different snapshot types.

An alternate approach  \cite{wash2013comm} applies POD to the state and sensitivity snapshots individually and concatenates the result to form the ROB.  If an orthogonal ROB is desired, a Gram-Schmidt-like procedure can be applied.  This approach has the advantages of: \begin{inparaenum}[(1)] \item avoiding heuristic snapshot weighting, \item individually controlling the contribution of each quantity to the ROB, and \item exploiting the optimality properties of POD for both the states and the sensitivities as both are required in the optimization context \end{inparaenum}.  Then, if $k_y^\text{state}$ is the dimension of the subspace spanned by the state vector basis and $k_y^\text{sens}$ is the dimension of the sensitivity subspace, the dimension of the ROB is $k_y = k_y^\text{state} + k_y^\text{sens}$.

\begin{algorithm}[htbp]
\caption{State-Sensitivity Proper Orthogonal Decomposition}
\begin{algorithmic}[1]\label{alg:state-sens-pod}
\REQUIRE State snapshot matrix, $\Xbold_\text{state} \in \Rbb^{N \times k_s^\text{state}}$, sensitivity snapshot matrix,  $\Xbold_\text{sens} \in \Rbb^{N \times k_s^\text{sens}}$, reduced state size, $k_y^\text{state}$ and reduced sensitivity size, $k_y^\text{sens}$
\ENSURE  ROB, $\Phibold$
  \STATE  Apply POD (Algorithm~\ref{alg:pod}) to $\Xbold_\text{state}$: $\text{POD}(\Xbold_\text{state},k_y^\text{state}) \rightarrow \Phibold_\text{state}$
  \STATE  Apply POD (Algorithm~\ref{alg:pod}) to $\Xbold_\text{sens}$: $\text{POD}(\Xbold_\text{sens},k_y^\text{sens}) \rightarrow \Phibold_\text{sens}$
  \STATE $\Phibold = \begin{bmatrix} \Phibold_\text{state} & \Phibold_\text{sens}\end{bmatrix}$
  \STATE (optional) Orthogonalize $\Phibold$ via QR (modified Gram-Schmidt)
\end{algorithmic}
\end{algorithm}

In \cite{arian2000trust}, the authors recognized the importance for optimization of having the ROM objective gradient 
approximate well the HDM objective gradient, particularly when aiming for convergence to the optimal solution.  
It was further stated that the incorporation of sensitivities in the POD basis without compromising the state 
approximation quality is an open question. An equally important issue is that of how state information in the ROB 
affects the approximation of a sensitivity. 

More recently, \cite{hay2009local} introduced two approaches for 
incorporating sensitivities into a ROB, one of which successfully improved the state approximation at the cost of 
increasing the basis size.  Additionally, a convergence theory for approximating sensitivities was provided in 
\cite{eftang2013approximation} in the context of the Empirical Interpolation Method (EIM). This theory implies
that sensitivity snapshots are unnecessary as the sensitivity approximation error will be driven to zero as the 
state approximation goes to zero.

In this work, the questions regarding the degradation of the state or sensitivity approximation by the
incorporation of extraneous snapshots are addressed by using \emph{minimum-residual} formulations for both the state 
and sensitivity equations.  In particular, the minimum-residual LSPG framework presented in Section~\ref{sec:rom-theory} 
ensures that appending sensitivity information to the ROB cannot degrade the state approximation, assuming a monotonic 
relationship between residual norm and error.  Additionally, a novel minimum-error framework for computing reduced 
sensitivities which ensures that appending state information to a ROB does not adversely affect the approximation of the 
sensitivities is presented in Section~\ref{sec:rom-sens}.


\section{OPTIMIZATION USING A PROGRESSIVE ROM} \label{sec:progromopt}
In this section, a framework for approximating the solution of the PDE-constrained optimization (\ref{opt:hdm}) using progressively constructed ROMs in a nonlinear trust region framework is presented.  The PDE-constrained optimization (\ref{opt:hdm}) is reduced using the ROM (\ref{eqn:rom}) as a surrogate for the high-dimensional PDE to produce the reduced optimization problem, which will be used interchangeably with the term \emph{ROM-constrained optimization},
\begin{equation} \label{opt:rom}
\begin{aligned}
& \underset{\ybold \in \Rbb^{k_y},~\mubold \in \boldsymbol{\Dcal}}{\text{minimize}}
& & \Jcal\left(\bar\wbold + \Phibold\ybold, \mubold\right) \\
& \text{subject to}
& & \cbold(\bar\wbold + \Phibold\ybold, \mubold) \leq 0 \\
& & &  \Psibold^T\Rbold(\bar\wbold + \Phibold\ybold, \mubold) = 0.
\end{aligned}
\end{equation} 
In the remainder of this section, the gradients of the ROM-constrained optimization problem are derived for Galerkin and LSPG ROMs, including a minimum-error formulation for the sensitivity computation.  Subsequently, the proposed progressive framework to approximate the solution of (\ref{opt:hdm}) is presented with all relevant algorithmic details made explicit.

At this point, the reader is reminded that the main objective of this paper is to demonstrate the accuracy of the proposed ROB-based computational technology.  Its potential for delivering the desired speedup will be demonstrated in future work, after it is equipped with hyperreduction.  Therefore, in the remainder, large dimensional quantities will appear in the computation of reduced 
quantities. 
 
\subsection{Sensitivities} \label{sec:rom-sens}
In the NAND framework, the reduced sensitivities $\displaystyle{\pder{\ybold}{\mubold}}$ can be derived by leveraging the fact that $\Psibold^T\Rbold(\bar\wbold + \Phibold\ybold(\mubold), \mubold) = 0$ is maintained throughout the solution process.  Employing a procedure identical to that used in Section~\ref{sec:nand}, in the general case where $\Psibold = \Psibold(\wbold_r(\ybold),\mubold)$, the reduced sensitivities are
\begin{equation} \label{eqn:rom-sens-gen}
\pder{\ybold}{\mubold} = -\left[\sum_{j=1}^N \Rbold_j \pder{\left(\Psibold^T\ebold_j\right)}{\wbold}\Phibold + \Psibold^T\pder{\Rbold}{\wbold}\Phibold\right]^{-1}\left[\sum_{j=1}^N \Rbold_j \pder{\left(\Psibold^T\ebold_j\right)}{\mubold} + \Psibold^T\pder{\Rbold}{\mubold}\right]
\end{equation}
where $\ebold_j$ is the $j$th canonical unit vector and the relation $\displaystyle{\pder{\wbold_r}{\ybold} = \Phibold}$ is used from (\ref{eqn:rom-assump}).  In the case of a Galerkin ROM, the left basis is state- and parameter-independent and the reduced sensitivities become
\begin{equation}\label{eqn:rom-sens-gal}
\pder{\ybold}{\mubold} = -\left[\Phibold^T\pder{\Rbold}{\wbold}\Phibold\right]^{-1}\left[\Phibold^T\pder{\Rbold}{\mubold}\right].
\end{equation}
The reduced sensitivities in the LSPG case are
\begin{equation}\label{eqn:rom-sens-lspg}
\pder{\ybold}{\mubold} = -\left[\sum_{j=1}^N \Rbold_j \Phibold^T\frac{\partial^2 \Rbold_j}{\partial \wbold \partial \wbold}\Phibold + \left(\pder{\Rbold}{\wbold}\Phibold\right)^T\left(\pder{\Rbold}{\wbold}\Phibold\right)\right]^{-1}\left[\sum_{j=1}^N \Rbold_j \Phibold^T\frac{\partial^2 \Rbold_j}{\partial \wbold \partial \mubold} + \left(\pder{\Rbold}{\wbold}\Phibold\right)^T\pder{\Rbold}{\mubold}\right]
\end{equation}
and require the (reduced) second-order sensitivities $\displaystyle{\Phibold^T\frac{\partial^2\Rbold_j}{\partial\wbold \partial \wbold}\Phibold}$ and $\displaystyle{\Phibold^T\frac{\partial^2 \Rbold_j}{\partial \wbold \partial \mubold}}$, which may not be available in black-box PDE solvers and expensive to compute with finite differences.  The higher-order terms may be neglected with minor loss in accuracy provided $\norm{\Rbold}$ is small, which will not, in general, be true in the context of parametric model order reduction.

In the case of hyperreduction, equations (\ref{eqn:rom-sens-gen}) - (\ref{eqn:rom-sens-lspg}) hold with the nonlinear terms replaced with the \emph{masked} nonlinear terms and the state vector $\wbold$ replaced with the restriction of the state vector to the sample mesh.  This ensures that all operations involved in the sensitivity computations are independent of the large dimension $N$.

Instead of computing the reduced sensitivities directly using one of (\ref{eqn:rom-sens-gen}) - (\ref{eqn:rom-sens-lspg}), approximate sensitivities can be constructed that minimize the distance to the corresponding HDM sensitivity in some norm $\Thetabold \succ 0$.  In this spirit, for each parameter, the desired sensitivity approximation is given by
\begin{equation} \label{eqn:red-sens-prob}
\widehat{\pder{\ybold}{\mubold_j}} = \underset{ \abold \in \Rbb^{k_y}}{\arg \min} ~\frac{1}{2} \norm{\pder{\wbold}{\mubold_j} - \Phibold\abold}_{\Thetabold}^2,
\end{equation}
where $\displaystyle{\widehat{\pder{\ybold}{\mubold}}}$ is to be used as a surrogate for $\displaystyle{\pder{\ybold}{\mubold}}$.
Using the expression for the HDM sensitivity (\ref{eqn:hdm-sens}) in (\ref{eqn:red-sens-prob}), the solution of the linear least-squares problems result in 
\begin{equation}\label{eqn:opt-red-sens1}
  \widehat{\pder{\ybold}{\mubold}} = -\left(\Thetabold^{1/2}\Phibold\right)^\dagger\Thetabold^{1/2}\pder{\Rbold}{\wbold}^{-1}\pder{\Rbold}{\mubold},
\end{equation}
where $\dagger$ denotes the Moore-Penrose pseudoinverse.  Selecting $\displaystyle{\Thetabold^{1/2} = \pder{\Rbold}{\wbold}}$, 
reduces (\ref{eqn:opt-red-sens1}) to
\begin{equation}\label{eqn:opt-red-sens}
  \widehat{\pder{\ybold}{\mubold}} = -\left(\pder{\Rbold}{\wbold}\Phibold\right)^\dagger\pder{\Rbold}{\mubold},
\end{equation}
which is equivalent to the optimality condition of
\begin{equation} \label{eqn:equiv-red-sens-prob}
  \underset{\abold \in \Rbb^{k_y}}{\text{minimize}} \quad \frac{1}{2} \norm{\pder{\Rbold}{\mubold_j} + \pder{\Rbold}{\wbold}\Phibold\abold}_2^2.
\end{equation}
This shows that choosing $\displaystyle{\Thetabold = \pder{\Rbold}{\wbold}^T\pder{\Rbold}{\wbold}}$ has the advantage of minimizing the state sensitivity error in the reduced subspace with respect to the $\Thetabold$-norm \emph{and} minimizing $\displaystyle{\norm{\frac{d\Rbold}{d\mubold}}}$ (used in derivation of HDM state sensitivity) in the span of $\displaystyle{\pder{\Rbold}{\wbold}\Phibold}$ with respect to the $2-$norm.  The corresponding normal equations are
\begin{equation}\label{eqn:opt-red-sens-normal-eqn}
  \left(\pder{\Rbold}{\wbold}\Phibold\right)^T\left(\pder{\Rbold}{\wbold}\Phibold\right)\widehat{\pder{\ybold}{\mubold}} = -\left(\pder{\Rbold}{\wbold}\Phibold\right)^T\pder{\Rbold}{\mubold}.
\end{equation}
These are equivalent to the LSPG sensitivity equations in (\ref{eqn:rom-sens-lspg}) with the second-order terms dropped.  In the context of optimization, the danger in using the minimum-error sensitivity computation is that the computed gradients will not be consistent with the functionals to which they correspond, which may result in convergence issues for the reduced optimization problem (\ref{opt:rom}).  The LSPG approach possesses the minimum-error sensitivity property (by dropping second-order terms), which approaches the true reduced sensitivity in regions of the parameter space where the ROM has high accuracy ($||\Rbold||$ small), i.e. near HDM samples.   

The minimum-error approach to computing reduced sensitivities ensures additional, unnecessary information in the ROB cannot degrade the sensitivity approximation, closing the loop in the discussion from Section~\ref{sec:offline}.  Additionally, using the state-sensitivity variant of POD in Algorithm~\ref{alg:state-sens-pod} and skipping the truncation step of the \emph{sensitivity} basis, we have $\displaystyle{\Phibold\widehat{\pder{\ybold}{\mubold}}(\wbold(\mubold),\mubold) = \pder{\wbold}{\mubold}(\wbold(\mubold),\mubold)}$ for the values of the parameters $\mubold$ at which the HDM was sampled, i.e. the reconstructed reduced senitivities exactly match the HDM sensitivities at training points.  This property follows directly from the minimum-error sensitivity formulation of the reduced sensitivities and the fact that the desired sensitivity is contained in the span of the ROB since the truncation step is skipped.


\subsection{ROM-constrained optimization with progressive reduced-basis construction}
The difficulty of training a static, global ROB which is expected to perform well throughout the parameter space, even for a small number of parameters, suggests the need for an adaptive procedure that merges the offline and online phases.  In the case of optimization, adaptive basis construction results in a ROM specialized for a region of the parameter space about the optimization trajectory.  In this work, the adaptive optimization procedure is driven by (\ref{opt:rom}) equipped with a residual-based nonlinear trust region
\begin{equation} \label{opt:rom-tr}
\begin{aligned}
& \underset{\ybold \in \Rbb^{k_y},~\mubold \in \boldsymbol{\Dcal}}{\text{minimize}}
& & \Jcal\left(\bar\wbold + \Phibold\ybold, \mubold\right) \\
& \text{subject to}
& & \cbold(\bar\wbold + \Phibold\ybold, \mubold) \leq 0 \\
& & &  \Psibold^T\Rbold(\bar\wbold + \Phibold\ybold, \mubold) = 0 \\
& & & \frac{1}{2}\norm{\Rbold(\bar\wbold + \Phibold\ybold, \mubold)}_2^2 \leq \epsilon.
\end{aligned}
\end{equation}
The proposed approach to ROM-constrained optimization is defined as a sequence of reduced optimization problems of the form (\ref{opt:rom-tr}).  The ROB, $\Phibold$, used in the $j$th reduced optimization problem is constructed using HDM samples at the solution of a subset of all previous reduced optimization problems.  Assuming a monotonic relationship between residual norm and error, proper selection of $\epsilon$ implies the trajectory of the $j$th reduced optimization problem will remain in a region of the parameter space where the corresponding ROM is accurate.  If the solution of the $j$th instance of (\ref{opt:rom-tr}), denoted $\mubold_j^*$, does not satisfy the optimality conditions of the PDE-constrained optimization problem (\ref{opt:hdm}), the ROB can be improved by sampling the HDM at $\mubold_j^*$.  As there is no formal convergence proof to guarantee that $\mubold_j^*$ will converge to a point satisfying the optimality conditions of (\ref{opt:hdm}), a heuristic convergence criteria 
\begin{equation} \label{eqn:converge}
  \norm{\mubold_{j-1} - \mubold_j} \leq \delta \norm{\mubold_j}
\end{equation}
is used, similar to \cite{afanasiev2001adaptive}.

With the fundamental approach introduced, subsequent sections will discuss crucial details that endow the proposed method with desirable properties, making the method viable for practical applications.  Among these details are: initial parameter value used for the $j$th optimization subproblem; the reference vector, $\bar\wbold$, and initial state, $\wbold^{(0)}$, used in the nonlinear iteration to solve the ROM (\ref{eqn:rom}) at each parameter configuration; selection of the residual upper bound, $\epsilon$, for the $j$th reduced optimization problem; and an algorithm for efficiently updating the ROB with new snapshots and reference vector without re-computing the SVD of the entire snapshot matrix, which grows as samples of the HDM are added.   For future reference, define the $k$th iteration of the $j$th reduced optimization problem as $\mubold_j^{(k)}$ and the solution of the $j$th reduced optimization problem as $\mubold_j^*$.  For notational convenience, define $\mubold_{-1}^* = \mubold_0^{(0)}$.

This section concludes with three definitions that are used in the remainder of this paper.  First, define $\Scal_j^{\mubold}$ as the set of parameter values at which the HDM has been sampled after the first $j+1$ instances of (\ref{opt:rom-tr}) 
\begin{equation} \label{eqn:s-mu}
  \Scal_j^{\mubold} = \{\mubold_{-1}^*, \mubold_0^*, \mubold_1^*, \dots, \mubold_j^*\}.
\end{equation}
Secondly, define the set of HDM solutions for each parameter in $\Scal_j^{\mubold}$
\begin{equation} \label{eqn:s-w}
  \Scal_j^{\wbold} = \{\wbold(\mubold_{-1}^*), \wbold(\mubold_0^*), \wbold(\mubold_1^*), \dots, \wbold(\mubold_j^*)\}.
\end{equation}
It may be advantageous to add several HDM samples to the training set before the optimization process is
started in order to improve the quality of the initial ROB.  This can be accomplished by  extending the definitions of $\Scal_j^{\mubold}$ and $\Scal_j^{\wbold}$ to include such samples. Further, define the ratio of actual objective reduction to predicted reduction, mimicking the corresponding quantity from traditional trust-region theory
\begin{equation} \label{eqn:actual-predict-ratio}
  \rho(\mubold, \nubold) = \frac{\Jcal(\wbold(\nubold),\nubold) - \Jcal(\wbold(\mubold),\mubold)}{\Jcal(\bar\wbold(\nubold)+\Phibold\ybold(\nubold),\nubold) - \Jcal(\bar\wbold(\mubold) + \Phibold\ybold(\mubold),\mubold)}.
\end{equation}
Then, define $\rho_j = \rho(\mubold_{j-1}^*, \mubold_j^*)$. The dependence of the reference vector, $\bar\wbold$, on the parameter vector, $\mubold$, is piecewise constant, as discussed in Section~\ref{sec:init-guess-rom}, validating the implicit assumption 
$\displaystyle{\pder{\bar\wbold}{\mubold} = 0}$ in Sections~\ref{sec:rom-theory} and~\ref{sec:rom-sens}.
\subsection{Initial guess for optimization subproblem: parameter space} \label{sec:init-guess-opt}
In the proposed progressive framework, a natural selection of the initial guess for the $j$th reduced optimization problem, $\mubold_j^{(0)}$, is the solution of the $j-1$ problem, $\mubold_{j-1}^*$.  However, if $\epsilon$ is too large in (\ref{opt:rom-tr}), the optimization trajectory may proceed into regions of the parameter space where it lacks accuracy.  In this case, the reduced optimization problem can converge to a parameter that \emph{increases} the true objective function compared to the initial guess, i.e. a ``worse'' point in terms of objective function minimization.

In anticipation of the case where $\epsilon$ is set too large, the initial guess for the reduced optimization problems is
the parameter at which the HDM was sampled that achieves the lowest  value of the objective function, namely,
\begin{equation} \label{eqn:init-guess-opt}
  \mubold_{j+1}^{(0)} = \arg \min_{\mubold \in \Scal_j^{\mubold}} \Jcal(\wbold(\mubold),\mubold).
\end{equation}
If the $j$th reduced optimization problem reduces the true objective function, the starting point for the subsequent problem 
is $\mubold_j^*$, as expected.  Notice that this choice of $\mubold_{j+1}^{(0)}$ does not require additional PDE solves as $\wbold(\mubold)$ is computed at every HDM sample.

In the case where $\epsilon$ is set too large, the HDM is still sampled at the solution, $\mubold_j^*$, and the ROB updated.  Therefore, the reduced optimization problem is not wasted because it adds information to the basis regarding areas of the parameter space to avoid, ensuring subsequent reduced optimization problems will not return to this point.  The heuristic regarding sampling the  ``worse'' point, while starting subsequent problems from the best point is an attempt to reduce the likelihood that an aggressive choice of $\epsilon$ will adversely affect convergence.

  
\subsection{Initial guess for ROM solve and reference vector selection: state space} \label{sec:init-guess-rom}
The solution of each reduced optimization problem (\ref{opt:rom-tr}) requires the solution of the nonlinear ROM equation (\ref{eqn:rom}) at each parameter value visited in the optimization procedure.  Regardless of the method used to solve the system of nonlinear equations, a starting point that is close to the true solution will significantly improve convergence to the solution.  Therefore, in the $j$th reduced optimization problem, given a parameter configuration $\mubold_j^{(k)}$ at which the ROM is to be solved, define the initial guess as
\begin{equation} \label{eqn:init-guess-solve}
  \wbold^{(0)} = \arg \min_{\wbold \in \Scal_j^{\wbold}} \norm{\Rbold(\wbold,\mubold_j^{(k)})}.
\end{equation}
This choice of the initial guess leverages the HDM samples taken thus far to use the best \emph{available} starting point, in the sense of the residual norm, which is usually used to assess convergence in nonlinear solvers.  From the previous section, the starting \emph{parameter configuration} for the $j$th optimization subproblem is also chosen from the set of HDM samples.  Therefore, if the optimization trajectory does not deviate too far from the starting point (\ref{eqn:init-guess-opt}), the initial state (\ref{eqn:init-guess-solve}) is expected to be a good starting point for the nonlinear solve. 

The reference vector is set to the initial condition, $\bar\wbold = \wbold^{(0)}$.  Selecting the reference vector in this way is a sufficient condition for consistency of a LSPG ROM \cite{carlberg2011gnat,carlberg2013gnat}, in the sense that the exact solution is recovered at sampled parameter values when POD is applied \emph{without} basis truncation.  Furthermore, it was observed in \cite{wash2012locbaseup} that using the initial condition as the reference vector outperforms other possible choices for small reduced-order bases. Perhaps the most significant implication of this choice of reference vector is that the ROM can exactly represent any state in the training set, regardless of truncation.  To elaborate on this point, consider $\mubold_{j+1}^{(0)} = \mubold_k^*$ with $k \leq j$.  From (\ref{eqn:init-guess-solve}), $\bar\wbold = \wbold^{(0)} = \wbold(\mubold_k^*)$ the MOR assumption in (\ref{eqn:rom-assump}) becomes $\wbold_r = \wbold(\mubold_k^*) + \Phibold\ybold.$  Taking $\ybold = 0$, the solution of $\Rbold(\wbold, \mubold_k^*) = 0$ is obtained, implying that the solution 
$\wbold(\mubold_k^*)$ is \emph{reachable} by the ROM independently of $\Phibold$.  The minimum-residual ROM formulation guarantees 
that the solution will be found, ignoring issues related to convergence.  The implication of this discussion is that the choices 
in (\ref{eqn:init-guess-opt}) - (\ref{eqn:init-guess-solve}) along with the affine representation in (\ref{eqn:rom-assump}) 
guarantee that the ROM is exact at the initial guess of each reduced optimization problem.  Therefore, the reduced optimization 
functionals precisely match those of the HDM, and the minimum-error reduced sensitivity framework guarantees that gradients 
will also match, provided that the sensitivity basis is not truncated in Algorithm~\ref{alg:state-sens-pod}.

Notice that modifying the reference vector changes the ROB due to the definition of state snapshots in Section~\ref{sec:offline}.  
Since the reference vector can potentially change at every iteration of each reduced optimization problem, it is 
paramount to avoid recomputing the ROB via POD with each new reference state.  Therefore, efficient, rank-1 basis updates are used 
following the work in \cite{wash2012locbaseup} to obtain a ROB corresponding to the \emph{translated} snapshots (see Algorithm~\ref{alg:svd-translate} and Section~\ref{sec:basis-updates} for additional details).

Using basis updates to ensure that the ROB corresponds (approximately) to the POD of the snapshots referenced by the 
appropriate initial condition has been shown to increase accuracy of a given ROM \cite{wash2012locbaseup}.  However, in the 
context of optimization, the modification to the ROB corresponds to changing the definition of the objective function and 
constraints \emph{within} the iterations of an optimization solve (\ref{opt:rom-tr}). Namely, the definition of the optimization 
functions are changing throughout the optimization process, which will likely have adverse effects on convergence.  
In practice, this does not show up as an issue because of the following observation.  At the $j$th reduced optimization problem, 
if $\mubold_j^{(0)}$ is close to $\mubold^*$, the solution of (\ref{opt:hdm}), it is unlikely that $\bar\wbold$ will change during 
the optimization procedure as $\mubold$ will not move far from $\mubold_j^{(0)}$: therefore, the ROM initial condition and 
reference vector will not move from $\wbold(\mubold_j^{(0)})$.  If $\mubold_j^{(0)}$ is far from $\mubold^*$, converging the 
optimization subproblem is not important provided progress is made toward $\mubold^*$; usual practice is to set an 
upper bound on the number of iterations allowed during the reduced optimization problems. Therefore, the potential issue 
with convergence of the reduced optimization problem is benign.

  
\subsection{Adaptive selection of residual upper bound} \label{sec:res-bound}
The selection of the residual upper bound, $\epsilon$, for a given reduced optimization problem, is discussed in this section.  As it is difficult to find a single value of $\epsilon$ a-priori that can be used for all reduced optimization problems, an adaptive approach is considered that modifies $\epsilon$ between successive reduced optimization solves.   Consider a value of the residual upper bound used in the $j$th reduced optimization problem, $\epsilon_j$; the value of $\epsilon_{j+1}$ is determined based on the value of $\rho_j$ (\ref{eqn:actual-predict-ratio}), the ratio of actual reduction of the objective function to predicted reduction.  It is clear that $\rho_j \approx 1$ indicates the ROM was a good model for the HDM in the prescribed trust region, while $\rho_j$ far from unity implies a poor model.  Therefore, the following update strategy from standard trust-region theory \cite{nocedal2006numerical} is proposed

\begin{equation} \label{eqn:eps-adapt}
	\epsilon' = \begin{cases} \displaystyle{\frac{1}{\tau}\epsilon} & \displaystyle{\rho_k \in \left[\frac{1}{2}, 2\right]}\\ \epsilon & \rho_k \in \left[\frac{1}{4}, \frac{1}{2}\right) \cup \left(2,4\right]\\ \tau\epsilon & \text{otherwise,} \end{cases}
\end{equation}
where $0 < \tau < 1$.  Solutions of the reduced optimization problem (\ref{opt:rom-tr}) that \emph{increase} the 
true objective function ($\rho_k < 0$) result in a reduction of the residual upper bound, $\epsilon$, by a factor of $\tau$, 
add an HDM sample to the training set, and ``reject'' the step in the sense that subsequent reduced optimization problems will be 
started from samples with lower values of the objective function as per (\ref{eqn:init-guess-opt}).
This is one among many possible adaptation strategies.  It is important to note that the adaptive $\epsilon$ only 
changes \emph{between} reduced optimization problems, and not within iterations of a given reduced problem.  
This ensures that the constraint definitions will not change between iterations in an optimization solve, unless $\bar\wbold$ changes as discussed in Section~\ref{sec:init-guess-rom}.

  The adaptive approach requires an initial value, $\epsilon_0$ to be manually defined, which can be done using problem-specific information, if available.  The goal of the adaptive $\epsilon$ procedure and the safe guards in Sections~\ref{sec:init-guess-opt}-\ref{sec:init-guess-rom}, is to make the progressive optimization process robust with respect to the initial value $\epsilon_0$.  If $\epsilon_0$ is smaller than necessary, the initial ROM will be highly accurate within the nonlinear trust region and result in a value of $\rho_1$ close to unity, which will cause the adaptive algorithm to increase the value of $\epsilon$ for subsequent reduced optimization problems.  If $\epsilon_0$ is too large, the optimization trajectory may venture into regions of the parameter space that cannot be captured by the ROM, resulting in a value of $\rho_1$ far from unity, causing the value of $\epsilon$ to be decreased for the next reduced optimization problem.  Notice that either under- or over-shooting the initial value of epsilon should only result in a few more HDM samples than would otherwise be necessary.


\subsection{Online algorithm for updating reduced-basis with new snapshots and reference vectors} \label{sec:basis-updates}
With each new HDM sample, a new ROB needs to be computed from the concatenation of the previous snapshot matrix with the new snapshots.  As more HDM samples are taken, the number of columns in the snapshot matrix grows by the number of snapshots collected per HDM sample.  For large problems, this quickly becomes prohibitively expensive as the cost of computing the SVD of a matrix scales quadratically with the number of columns \cite{golub2012matrix}.  Additionally, new reference vectors computed in each reduced 
optimization problem from (\ref{eqn:init-guess-solve}) require computation of a new ROB with the state snapshots translated by 
the new $\bar\wbold$.  This will also introduce a significant bottleneck as (\ref{eqn:init-guess-solve}) allows for the possibility 
of changing $\bar\wbold$ at each iteration of each reduced optimization problem, although in practice it is not this often.
In this section, an extension of the work described in \cite{wash2012locbaseup} that applies low-rank SVD updates \cite{brand2006fast} is presented to economically update a ROB with new snapshots and a new reference vector by leveraging the SVD of previous snapshot matrices.

The low-rank SVD update algorithm in \cite{brand2006fast} is used to cheaply compute the thin SVD of $\Xbold + \Abold\Bbold^T$ from the SVD of $\Xbold$.  With proper selection of $\Xbold$, $\Abold,$ and $\Bbold$, this algorithm can be tailored to the cases of appending vectors to the data matrix and translating the columns \cite{brand2006fast}, see Algorithms~\ref{alg:svd-append} and~\ref{alg:svd-translate}, respectively, in the Appendix.

In the context of the state-sensitivity POD (Algorithm~\ref{alg:state-sens-pod}), a new reference vector and additional state and sensitivity snapshots imply that the state basis must incorporate the additional state snapshots and reference vector and the sensitivity basis must incorporate the additional sensitivity snapshots.  If $\Xbold_\text{state}$ is the original state snapshot matrix and $\bar\wbold$ is the original reference vector, then $\Phibold_\text{state} = \mathrm{POD}(\Xbold_\text{state} - \bar\wbold\onebold^T)$, where $\onebold$ is the vector of ones of appropriate size.  Similarly, the sensitivity basis is $\Phibold_\text{sens} = \mathrm{POD}(\Xbold_\text{sens})$, where $\Xbold_\text{sens}$ is the sensitivity snapshot matrix.  With a new reference vector, $\widehat\wbold$, and snapshot matrices, $\Ybold_\text{state}$, and $\Ybold_\text{sens}$, bases of the form $\widehat\Phibold_\text{state} = \mathrm{POD}\left(\begin{bmatrix} \Xbold_\text{state} & \Ybold_\text{state} \end{bmatrix} - \widehat\wbold\onebold^T\right)$ and  $\widehat\Phibold_\text{sens} = \mathrm{POD}\left(\begin{bmatrix} \Xbold_\text{sens} & \Ybold_\text{sens} \end{bmatrix}\right)$ are desired.  This can be achieved without directly recomputing the POD of the modified snapshot matrices by applying low-rank basis updates.  The desired bases can be written as
\begin{align*}
\widehat\Phibold_\text{state} &= \mathrm{POD}\left(\begin{bmatrix} \Xbold_\text{state} - \bar\wbold\onebold^T + \left(\bar\wbold - \widehat\wbold\right)\onebold^T & \Ybold_\text{state} - \widehat\wbold\onebold^T \end{bmatrix}\right) \\
\widehat\Phibold_\text{sens} &= \mathrm{POD}\left(\begin{bmatrix} \Xbold_\text{sens} & \Ybold_\text{sens} \end{bmatrix}\right).
\end{align*}
This form makes it clear that a \emph{translation} update with the vector $\bar\wbold - \widehat\wbold$ can be applied to $\Phibold_\text{state}$ and an \emph{append} update with the vectors $\Ybold_\text{state} - \widehat\wbold\onebold^T$ to yield the desired $\widehat\Phibold_\text{state}$.  Similarly, $\widehat\Phibold_\text{sens}$ is obtained by applying an \emph{append} update with the matrix $\Ybold_\text{sens}$.  This procedure is summarized in Algorithm~\ref{alg:svd-append-trans}.  To obtain the updated ROB with the new reference vector $\widehat\wbold$ without additional snapshots, as required in Section~\ref{sec:init-guess-rom}, Algorithm~\ref{alg:svd-append-trans} can be used with $\Ybold_\text{state} = \Ybold_\text{sens} = \emptyset$.

\begin{algorithm}
  \caption{Updating ROB from state-sensitivity POD: appending vectors and changing reference vector} \label{alg:svd-append-trans}
  \begin{algorithmic}[1]
    \REQUIRE State and sensitivity bases: $\Phibold_\text{state}$, $\Phibold_\text{sens}$ as well as corresponding singular values ($\Sigmabold_\text{state}$ and $\Sigmabold_\text{sens}$) and right singular vectors ($\Vbold_\text{state}$ and $\Vbold_\text{sens}$); original reference vector, $\bar\wbold$; new reference vector, $\widehat{\bar\wbold}$; new state snapshots, $\Ybold_\text{state}$; new sensitivity snapshots, $\Ybold_\text{sens}$
    \ENSURE Updated ROB, $\widehat\Phibold $
    \STATE Update $\Phibold_\text{sens}$ with additional snapshots, $\Ybold_\text{sens}$ via Algorithm~\ref{alg:svd-append} $\rightarrow \widehat\Phibold_\text{sens}$
    \STATE Update $\Phibold_\text{state}$ with new reference vector via Algorithm~\ref{alg:svd-translate} applied to $\bar\wbold - \widehat{\bar\wbold} \rightarrow \bar\Phibold_\text{state}$
    \STATE Update $\bar\Phibold_\text{state}$ with additional snapshots, $\Ybold_\text{state} - \widehat{\bar\wbold}\onebold^T \rightarrow \widehat\Phibold_\text{state}$
    \STATE Define $\widehat\Phibold = \begin{bmatrix} \widehat\Phibold_\text{state} & \widehat\Phibold_\text{sens}\end{bmatrix}$
    \STATE (optional) Orthogonalize $\widehat\Phibold$ via QR (modified Gram-Schmidt)
  \end{algorithmic}
\end{algorithm}

It is important to mention that the low-rank SVD updates in \cite{brand2006fast} assume the entire economy SVD is available to perform the update.  In our application, only the truncated SVD may be available, resulting in an error in the update.  It can be shown that, in the case of column translation, the error is proportional to the largest neglected singular value (usually small in MOR applications).  Additionally, it can be shown that the error in the update does not noticeably affect the quality of the basis.  Finally, it can also be shown that the thin SVD update of Brand can be applied in the context of hyperreduction, where only \emph{masked} singular factors are available.


\subsection{Summary: proposed ROM-constrained optimization algorithm}
Section~\ref{sec:progromopt} concludes with a summary of the proposed approach to ROM-constrained optimization using progressively-constructed ROMs, including all relevant implementation details.

\begin{algorithm}[ht]
  \caption{PDE-constrained optimization using progressively-constructed ROMs} \label{alg:progromopt}
  \begin{algorithmic}[1]
    \REQUIRE Initial guess for optimization, $\mubold_0^{(0)}$, initial residual upper bound, $\epsilon_0$, convergence tolerance, $\delta$ and maximum number of reduced optimization solves, $N_\text{max}$
    \ENSURE Approximation to solution of (\ref{opt:hdm}), $\mubold_r^*$ 
    \STATE Set $\Scal_{-1}^{\mubold} = \{\}$ and $\Scal_{-1}^{\wbold} = \{\}$
    \STATE Define $\mubold_{-1}^* = \mubold_0^{(0)}$
    \FOR{$j = 0, 1, \dots, N_\text{max}$} 
     \STATE Compute $\wbold(\mubold_{j-1}^*)$ and $\pder{\wbold}{\mubold}(\mubold_{j-1}^*)$ (sample HDM)
     \STATE Set $\Scal_{j} = \Scal_{j-1}^{\mubold} \cup \{\mubold_{j-1}^*\}$ and $\Scal_{j}^{\wbold} = \Scal_{j-1}^{\wbold} \cup \{\wbold(\mubold_{j-1}^*)\}$
     \IF{j == 0}
       \STATE Apply Algorithm~\ref{alg:state-sens-pod} to state and sensitivity snapshots, $\wbold(\mubold_{j-1}^*)$ and $\pder{\wbold}{\mubold}(\mubold_{j-1}^*) \rightarrow \Phibold$
       \ELSE
         \STATE Compute $\epsilon_j$ using (\ref{eqn:eps-adapt})
         \STATE Apply Algorithm~\ref{alg:svd-append-trans} to update the ROM, $\Phibold$ with additional state and sensitivity snapshots $\wbold(\mubold_{j-1}^*)$ and $\pder{\wbold}{\mubold}(\mubold_{j-1}^*) \rightarrow \Phibold$
     \ENDIF
     \STATE Compute $\mubold_j^{(0)}$ using (\ref{eqn:init-guess-opt})
      \STATE Solve (\ref{opt:rom-tr}) with ROB $\Phibold$ and residual upper bound $\epsilon_j$; initial condition and offset, $\bar\wbold$, determined via (\ref{eqn:init-guess-solve}) for each ROM solve (use Algorithm~\ref{alg:svd-append-trans} with $\Ybold_\text{state} = \Ybold_\text{sens} = \emptyset$ to update the ROB with the new $\bar\wbold$ when necessary) $\rightarrow \mubold_{j}^*$
       \IF{$\norm{\mubold_j^* - \mubold_{j-1}^*}_2 < \delta\norm{\mubold_j^*}_2$}
         \STATE $\mubold_r^* = \mubold_j^*$
         \STATE break
       \ENDIF
    \ENDFOR
  \end{algorithmic}
\end{algorithm}

\section{APPLICATIONS} \label{sec:app}

Here, the approach to PDE-constrained optimization using a progressively-constructed parametric ROM described in 
Section~\ref{sec:progromopt} is illustrated with the solution of an aerodynamic \emph{shape} optimization problem.  
Specifically, pressure distribution discrepancy in the subsonic regime is used to recover the shape of the RAE2822 airfoil
starting from that of the NACA0012 airfoil.

\subsection{Shape parametrization and problem setup}

A plethora of shape parametrization techniques exist \cite{samareh1999survey,anderson2012parametric}, each with strengths 
and weaknesses. They typically trade-off between efficiency and flexibility. A subset of these techniques has been studied 
in the context of model order reduction \cite{rozza2010model}. In this work, the SDESIGN software \cite{mautesb1,maute2006fem}, based on the design element approach \cite{imam1982three,farin1996curves}, is used for shape parameterization.  
In this approach, an initial shape is mapped onto a set of design (or control) elements with nodal displacement degrees of freedom.
These are used to morph the given shape into a different one by applying to it a displacement field represented by shape functions
generated by the design elements and their nodal degrees of freedom. 
Here, a single ``cubic'' design element is used to parametrize the 
NACA0012 airfoil using this approach. Such a design element has 8 control nodes. They are used to define cubic Lagrangian 
polynomials for describing the displacement fields along the horizontal edges of the element, and linear functions for 
representing the displacement fields along its vertical ones. For this application, the set of admissible shapes generated by 
this representation is restricted by constraining the control nodes to move in the \emph{vertical} direction only.  
This results in a parameterization with 8 variables where each of them represents the \emph{displacement} of 
a control node in the vertical direction. The case where all parameters are equal, $\mubold = c\onebold$ for $c \in \Rbb$, 
corresponds to a rigid translation in the vertical direction. Because such a translation does not affect the definition of 
a shape, it is eliminated by constraining one of the displacement variables to zero.  Furthermore, because
the control nodes are allowed to move only in the vertical direction, rigid rotations are automatically eliminated.

\begin{figure}[h]
  \centering
  \begin{subfigure} [b]{0.45\textwidth}
    \centering
    \includegraphics[width=\textwidth]{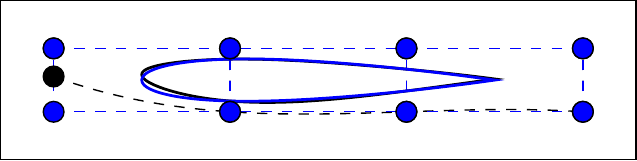}
    \caption{$\mubold(1)$ = 0.1}
  \end{subfigure}
    \begin{subfigure} [b]{0.45\textwidth}
    \centering
    \includegraphics[width=\textwidth]{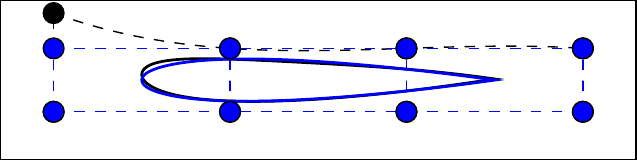}
    \caption{$\mubold(2)$ = 0.1}
  \end{subfigure}\vspace{0.5cm}\\
  \begin{subfigure} [b]{0.45\textwidth}
    \centering
    \includegraphics[width=\textwidth]{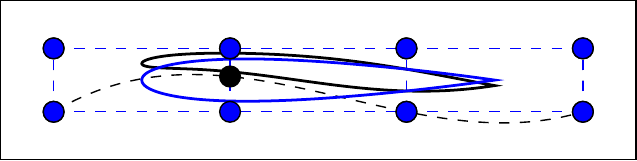}
    \caption{$\mubold(3)$ = 0.1}
  \end{subfigure}
    \begin{subfigure} [b]{0.45\textwidth}
    \centering
    \includegraphics[width=\textwidth]{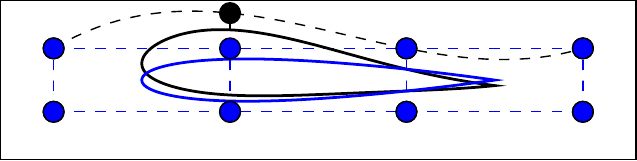}
    \caption{$\mubold(4)$ = 0.1}
  \end{subfigure} \vspace{0.5cm}\\
  \begin{subfigure} [b]{0.45\textwidth}
   \centering
    \includegraphics[width=\textwidth]{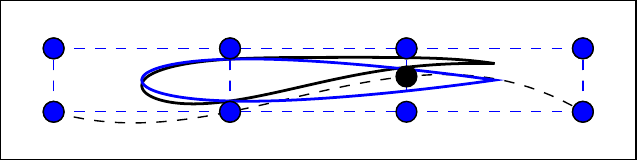}
    \caption{$\mubold(5)$ = 0.1}
  \end{subfigure}
    \begin{subfigure} [b]{0.45\textwidth}
    \centering
    \includegraphics[width=\textwidth]{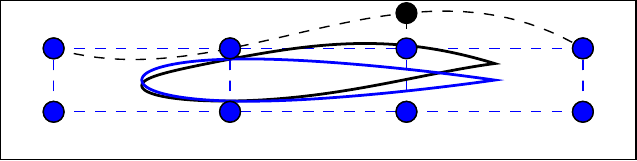}
    \caption{$\mubold(6)$ = 0.1}
  \end{subfigure}\vspace{0.5cm}\\
  \begin{subfigure} [b]{0.45\textwidth}
    \centering
    \includegraphics[width=\textwidth]{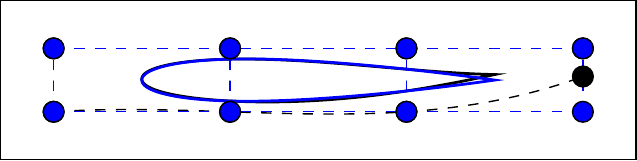}
    \caption{$\mubold(7)$ = 0.1}
  \end{subfigure}
  \begin{subfigure} [b]{0.45\textwidth}
    \centering
    \includegraphics[width=\textwidth]{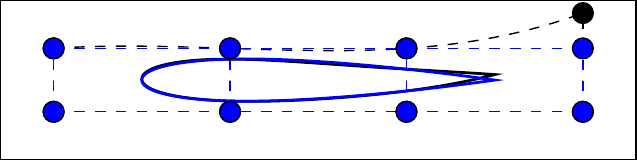}
    \caption{$\mubold(8)$ = 0.1}
  \end{subfigure}
  \caption{Shape parametrization of a NACA0012 airfoil using a \emph{cubic} design element (the notation $\mubold(i)$ designates the $i$-th component of the vector $\mubold$ which refers to the $i$-th displacement degree of freedom of the shape parameterization)}
  \label{fig:shape-param}
\end{figure}

A visualization of the vertices of the design element and the deformation induced by perturbing each design variable is
given in Figure~\ref{fig:shape-param}. While SDESIGN is used to deform the surface nodes of the airfoil, a robust mesh motion 
algorithm based on a structural analogy is used to deform the surrounding body-fitted CFD mesh accordingly.

\begin{figure}[ht]
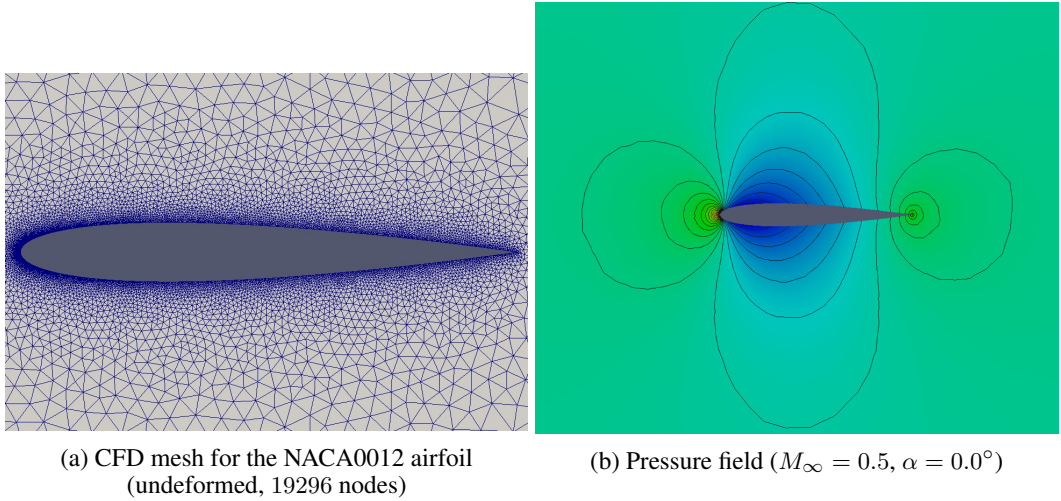

  \centering
  \begin{subfigure} [b]{0.48\textwidth}
    \centering
    \includegraphics[width=\textwidth]{results/exo/{{naca0012Medium.undef.mesh}}}
    \caption{CFD mesh for the NACA0012 airfoil \\ (undeformed, $19296$ nodes)} \label{fig:naca0012-undef-mesh}
   \end{subfigure}
   \begin{subfigure} [b]{0.48\textwidth}
     \centering
     \includegraphics[width=\textwidth]{results/exo/{{naca0012Medium.undef.M0p5.a0p0.b0p0.withcontours}}}
     \caption{Pressure field ($M_{\infty} = 0.5$, $\alpha = 0.0^\circ$) \\ \phantom{place holder}} \label{fig:naca0012pdistM0.5}
   \end{subfigure}
   \caption{NACA0012 mesh and pressure distribution at Mach 0.5 and zero angle of attack} \label{fig:naca0012}
\end{figure}

The flow over the airfoil is modeled using the compressible Euler equations, and these are solved numerically using
AERO-F \cite{AEROF}. Because this flow solver is three-dimensional, the two-dimensional fluid domain around the airfoil
is represented as a slice of a three-dimensional domain. This slice is discretized using a body-fitted CFD mesh with $54816$ 
tetrahedra and $19296$ nodes (Figure~\ref{fig:naca0012-undef-mesh}). Specifically, the flow equations are semi-discretized 
by AERO-F on this CFD mesh using a second-order finite volume method based on Roe's flux \cite{Roe}. 

For each airfoil configuration generated during the iterative optimization
procedure, the steady state solution of the flow problem is computed iteratively using pseudo-time-integration. 
For this purpose, each sought-after steady state solution is initialized using the \emph{best} previously computed steady state 
solution available in the database\footnote{In this context, the database refers to the flow solutions computed for all shapes 
previously visited by the optimization trajectory.}. The best steady state solution is defined here as that steady state solution 
available in the database which, for the given airfoil configuration, minimizes the residual of the discretized steady state 
Euler equations. Because the database of steady state flow solutions is initially empty, the iterative computation of the steady
state flow over the initial shape --- in this case, that of the NACA0012 airfoil --- is initialized with the uniform flow solution.


The framework for model order reduction described in Section~\ref{sec:rom-theory} is used to solve the aerodynamic shape optimization problem.  At each HDM sample, the steady state solution and sensitivities with respect to shape parameters are computed and used as snapshots.  As the chosen shape parameterization has 8 parameters, 9 snapshots are generated per HDM sample: one snapshot 
corresponding to the steady state solution, and 8 to the sensitivity of this solution with respect to all 8 shape parameters. 
A ROB is extracted from 
these snapshots using the state-sensitivity variant of the POD method described in Section~\ref{sec:offline}.  
Because very few snapshots are generated for this problem, the truncation step in the POD algorithm 
(Algorithm~\ref{alg:pod}) is skipped. Consequently, the size of the constructed ROB is $k_y = 9s$, where $s$ is the number 
of sampled HDMs.

Each ROM problem is solved using the Gauss-Newton method equipped with a backtracking linesearch algorithm
to minimize (\ref{eqn:min-res}). For this purpose, the initial condition and reference state vector are chosen as 
in (\ref{eqn:init-guess-solve}). The python interface to the SNOPT \cite{gill2002snopt} software, pyOpt \cite{pyopt-paper}, is 
used to solve the optimization problem itself.

\begin{figure}[h]
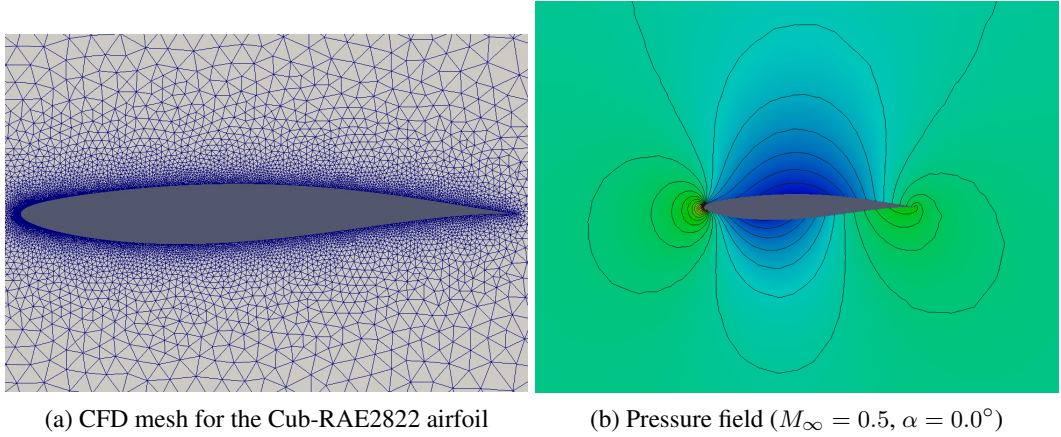

  \centering
  \begin{subfigure} [b]{0.48\textwidth}
    \centering
    \includegraphics[width=\textwidth]{results/exo/{{naca0012Medium.rae2822.mesh}}}
    \caption{CFD mesh for the Cub-RAE2822 airfoil} \label{fig:rae2822mesh}
   \end{subfigure}
   \begin{subfigure} [b]{0.48\textwidth}
     \centering
     \includegraphics[width=\textwidth]{results/exo/{{naca0012Medium.rae2822.M0p5.a0p0.b0p0.withcontours}}}
     \caption{Pressure field ($M_{\infty} = 0.5$, $\alpha = 0.0^\circ$)} \label{fig:rae2822pdist}
   \end{subfigure}
   \caption{Cub-RAE2822 mesh and pressure isolines computed at Mach 0.5 and zero angle of attack} \label{fig:rae2822}
\end{figure}

At this point, it is noted that since the exact profile of the RAE2822 airfoil does not lie in the space of admissible airfoil 
profiles defined by the cubic design element parametrization described above, it is approximated by the closest admissible profile.  
This approximation is referred to in the remainder of this section as the Cub-RAE2822 airfoil. It is graphically
depicted in Figure~\ref{fig:rae2822} which also shows the pressure isolines computed for this airfoil at the
free-stream Mach number $M_{\infty} = 0.5$ and angle of attack $\alpha = 0.0^\circ$.

\subsection{Subsonic inverse design} \label{sec:inv-design}
The free-stream conditions of interest are set to the subsonic Mach number $M_{\infty} = 0.5$ and zero angle of attack
($\alpha = 0^{\circ}$), and the following optimization problem is considered
\begin{equation} \label{opt:match}
\begin{aligned}
& \underset{\wbold \in \Rbb^N,~\mubold \in \boldsymbol{\Dcal}}{\text{minimize}}
& & \frac{1}{2}\norm{\pbold(\wbold(\mubold)) - \pbold(\wbold(\mubold^\text{RAE2822}))}_2^2 \\
& \text{subject to}
& & \mubold(3) = 0 \\
& & &  \boldsymbol{\ell} \leq \mubold \leq \ubold \\
& & & \Rbold(\wbold, \mubold) = 0
\end{aligned}
\end{equation} 
where $\pbold(\wbold)$ is the vector of nodal pressures, and $\mubold^\text{RAE2822}$ designates the 
parameter solution vector morphing the NACA0012 airfoil into the Cub-RAE2822 airfoil. 
The first constraint is introduced to eliminate the
rigid body translation in the vertical direction as discussed in the previous section.
The box constraints are introduced to prohibit the optimization trajectory from going through highly distorted shapes 
for which it is difficult to compute steady state flow solutions.

\begin{figure}[h]
   \centering
   \includegraphics[width=\textwidth]{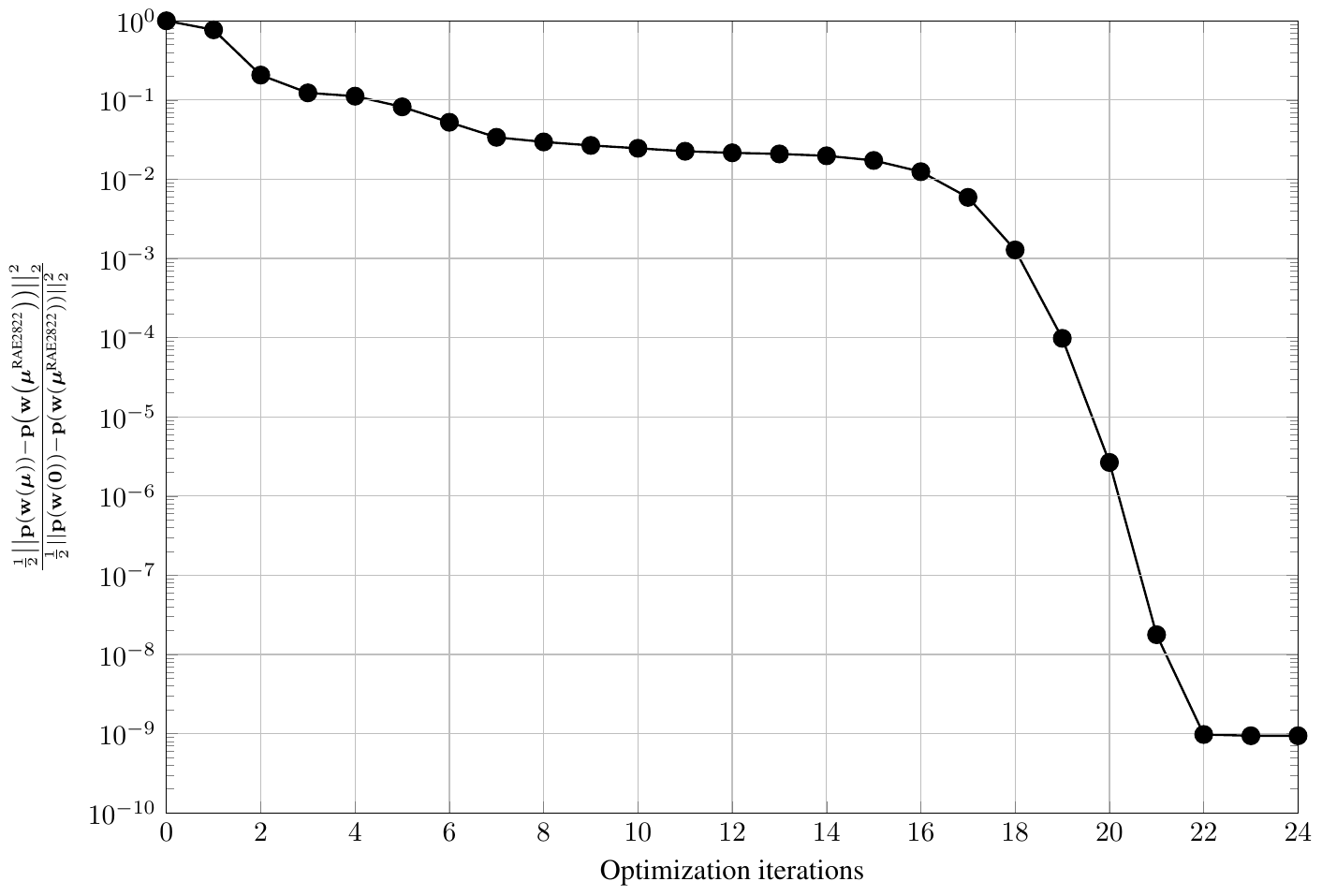}
   \caption{Progression of the objective function during the HDM-based optimization.  The initial guess is defined as the $0$th optimization iteration.}
   \label{fig:match-hdm-obj}
\end{figure}

\begin{figure}[h]
  \centering
   \includegraphics[width=\textwidth]{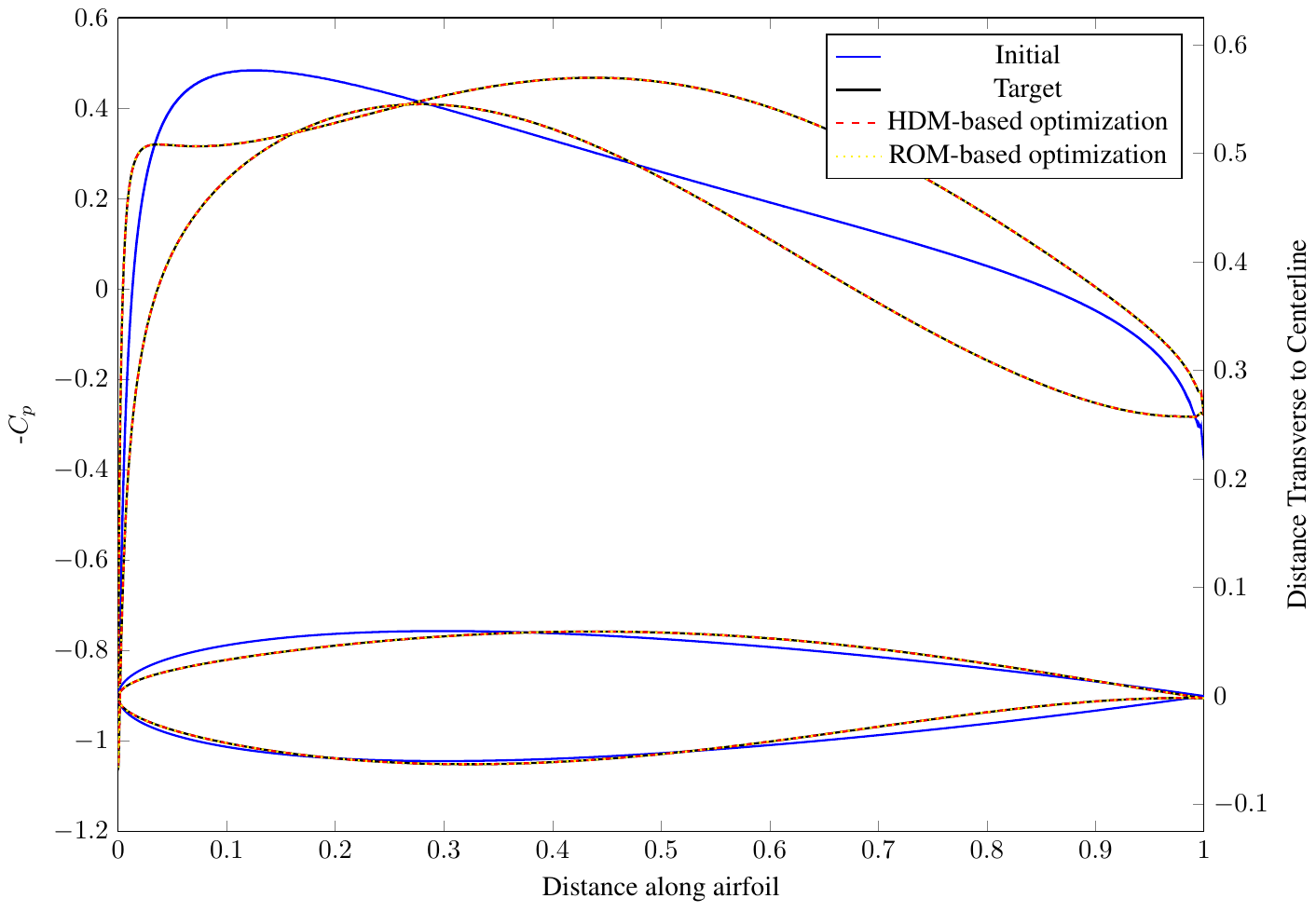}
   \caption{Subsonic inverse design of the airfoil Cub-RAE2822: initial shape (NACA0012) and associated $C_p$ function, and final 
	   shape (Cub-RAE2822) and associated $C_p$ functions delivered by the HDM- and ROM-based optimizations, respectively} \label{fig:match-cp-disp}
\end{figure}

To obtain a reference solution that can be used for assessing the performance of the proposed ROM-based optimization method,
problem (\ref{opt:match}) is first solved using the HDM as the constraining PDE. In this case, the optimizer is found to reduce 
the initial value of the objective function by 9 orders of magnitude, before numerical difficulties cause 
it to terminate (see Figure~\ref{fig:match-hdm-obj}). Relevant statistics associated with this HDM-based reference solution of the
optimization problem are gathered in Table~\ref{tab:summary}. Essentially, 24
optimization iterations are required for obtaining a solution $\mubold^\text{RAE2822}$ with a relative error well below $0.1\%$. 
These iterations incur a total of 29 HDM queries (including those associated with the linesearch iterations). 
Figure~\ref{fig:match-cp-disp} shows that the pressure coefficient function associated with this reference solution matches 
very well the target pressure coefficient function.   

Next, the framework for ROM-constrained optimization described in Section~\ref{sec:progromopt}, which is
summarized in Algorithm~\ref{alg:progromopt}, is applied to the solution of problem (\ref{opt:match-red}) 
via the solution of a sequence of reduced optimization problems of the form
\begin{equation} \label{opt:match-red}
\begin{aligned}
& \underset{\wbold \in \Rbb^N,~\mubold \in \boldsymbol{\Dcal}}{\text{minimize}}
& & \frac{1}{2}\norm{\pbold(\bar\wbold + \Phibold \ybold(\mubold)) - \pbold(\wbold(\mubold^\text{RAE2822}))}_2^2 \\
& \text{subject to}
& & \mubold(3) = 0 \\
& & &  \boldsymbol{\ell} \leq \mubold \leq \ubold \\
& & & \Psibold^T\Rbold(\bar\wbold+\Phibold\ybold(\mubold), \mubold) = 0 \\ 
& & & \frac{1}{2}\norm{\Rbold(\bar\wbold+\Phibold\ybold(\mubold), \mubold)}_2^2 \leq \epsilon.
\end{aligned}
\end{equation} 

The HDM is sampled at the initial configuration and the resulting 9 snapshots are used to build a ROB using the 
state-sensitivity POD method summarized in Algorithm~\ref{alg:state-sens-pod}, without truncation.  
The resulting ROB is used to solve (\ref{opt:match-red}).  Indeed, as the minimum-error sensitivity 
computation described in Section~\ref{sec:rom-sens} is not consistent with the true reduced sensitivities for large 
residuals, convergence of the optimization problem is not guaranteed. To address this issue, an upper bound is set 
on the number of optimization iterations (25 in this case) and the goal of the reduced optimization problem is set to finding an 
improvement to the current solution before updating the ROB.  The HDM is also sampled at the termination point of each 
reduced optimization problem yielding $9$ additional snapshots which are appended to the ROB using Algorithm~\ref{alg:svd-append-trans}.  Linear independence of the basis is maintained by truncating vectors with corresponding singular values below some tolerance.  For the present application, such truncation was not necessary as the snapshots added to the ROB at a given iteration were not contained in the span of the snapshots from previous iterations.  When solving for the flow using 
a given ROM, the initial guess and snapshot translation vector are defined in Section \ref{sec:init-guess-rom} as the best HDM sample.  Additionally, an aggressive choice of $\tau = 0.1$ was used for the nonlinear trust region radius bound adaptation in (\ref{eqn:eps-adapt}).

Using only 7 HDM samples, the progressive ROM optimization framework reduces the initial pressure discrepancy by 18 orders of 
magnitude, to essentially machine zero. Interestingly, this is 4 times fewer HDM queries than required by 
the HDM-based optimization. Figure~\ref{fig:match-cp-disp} shows that both the shape of the Cub-RAE2822 airfoil
and the associated pressure coefficient at the chosen free-stream conditions match the target shape and pressure coefficient
function very well.

\begin{figure}[ht]
   \centering
   \includegraphics[width=\textwidth]{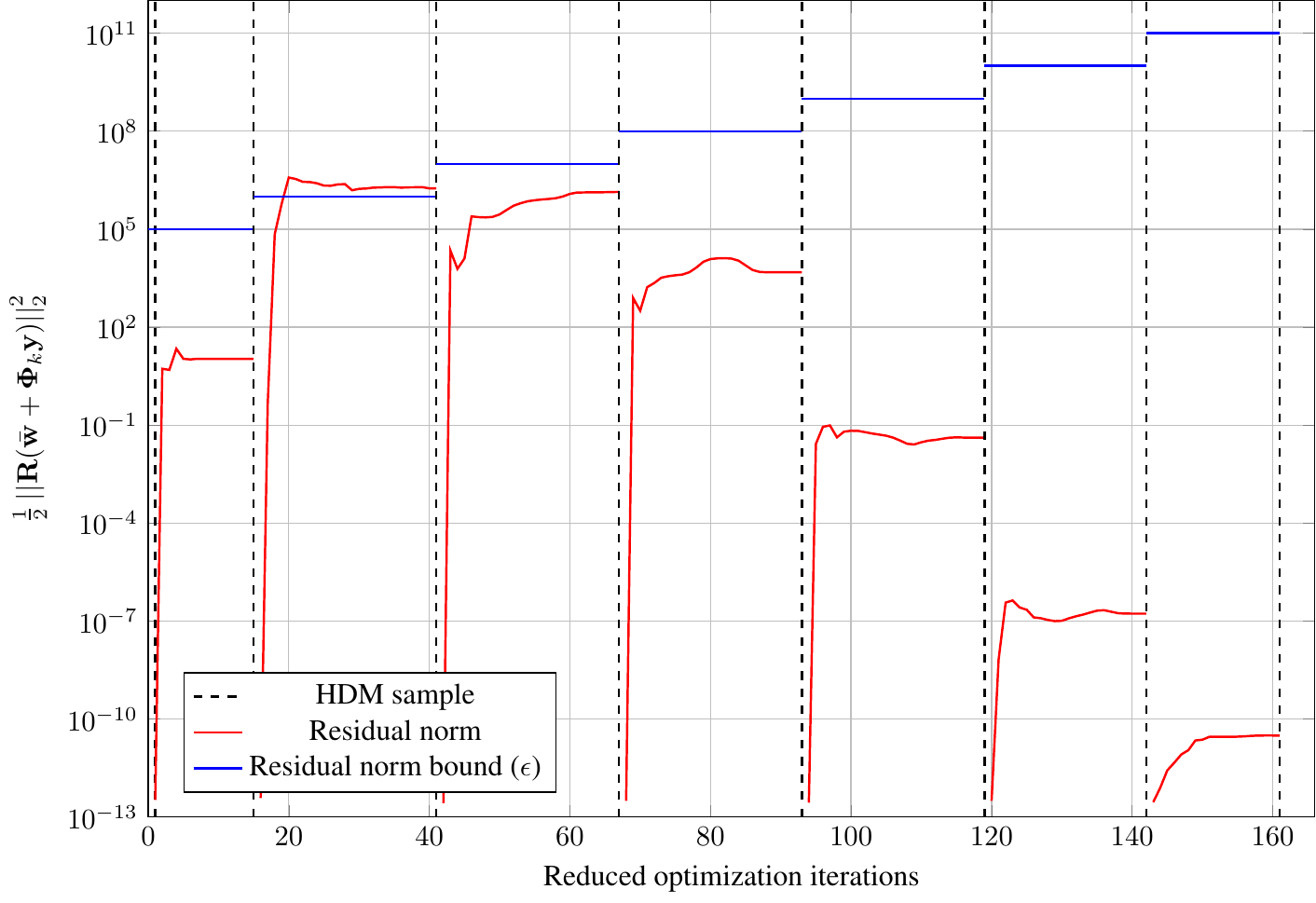}
   \caption{Progression of HDM residual: dashed line indicates an HDM sample and a subsequent update of the ROB}
   \label{fig:match-rom-fluxnorm}
\end{figure}

Surprisingly, the ROM-based optimization process achieves a \emph{lower} value of the objective function than the HDM-based one.  This can be traced to convergence tolerance on the HDM sensitivity analysis.  The HDM-based sensitivities are obtained by solving the multiple right-hand side linear system of equations in (\ref{eqn:hdm-sens}) using GMRES.  The convergence tolerance is $\norm{\Abold\xbold - \bbold}_2 \leq \gamma\norm{\bbold}_2$ for solving the linear system of equations $\Abold\xbold = \bbold$, with $\gamma = 10^{-10}$ in this case.  If $\norm{\bbold}$ is large ($\bbold = \displaystyle{\pder{\Rbold}{\mubold}}$ in this case), the convergence requirement may be rather flexible.  Conversely, the ROM-based sensitivities in (\ref{eqn:opt-red-sens}) are solved directly using the QR factorization, where no tolerances are involved.  Recall from Section~\ref{sec:rom-sens} that the minimum-error reduced sensitivities approach the true sensitivities for LSPG projection as the HDM residual approaches zero.  Figure~\ref{fig:match-rom-fluxnorm} verifies that the HDM residual is small after 6 HDM samples are taken, which implies the minimum-error ROM sensitivities are (nearly) consistent with the true ROM sensitivities.  This consistency will guarantee convergence of the reduced optimization problem when using a globally-convergent optimization solver.  Additionally, the small HDM residual implies that the ROM is highly accurate in this region, making it likely that the reduced optimization problem will converge to a point close to the true optimum.

\begin{figure}[h]
   \centering
   \includegraphics[width=\textwidth]{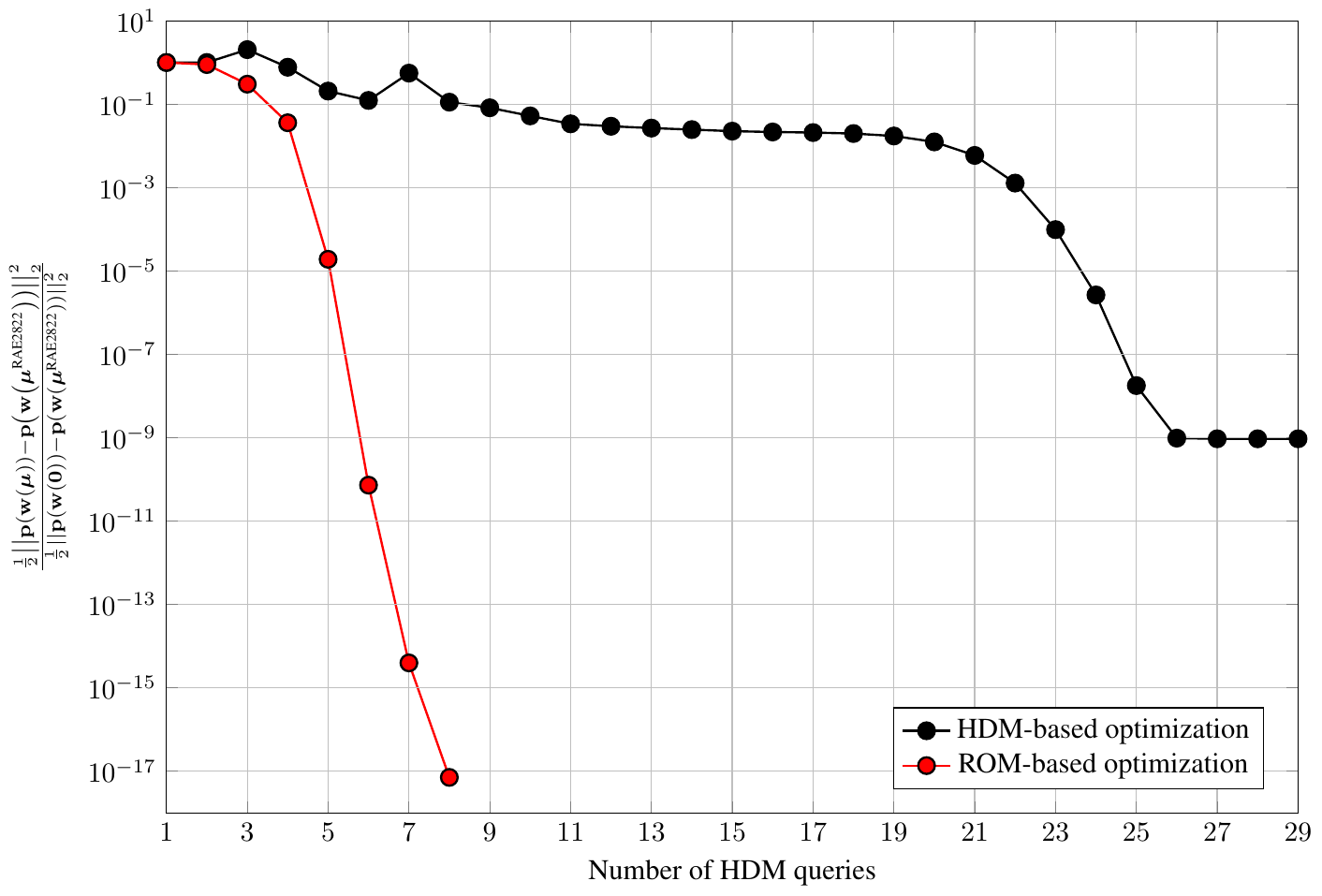}
   \caption{Objective function versus number of queries to the HDM: ROM-based optimization (red) and HDM-based optimization (black)}
   \label{fig:match-hdmrom-obj}
\end{figure}

\begin{figure}[h]
     \centering
     \includegraphics[width=\textwidth]{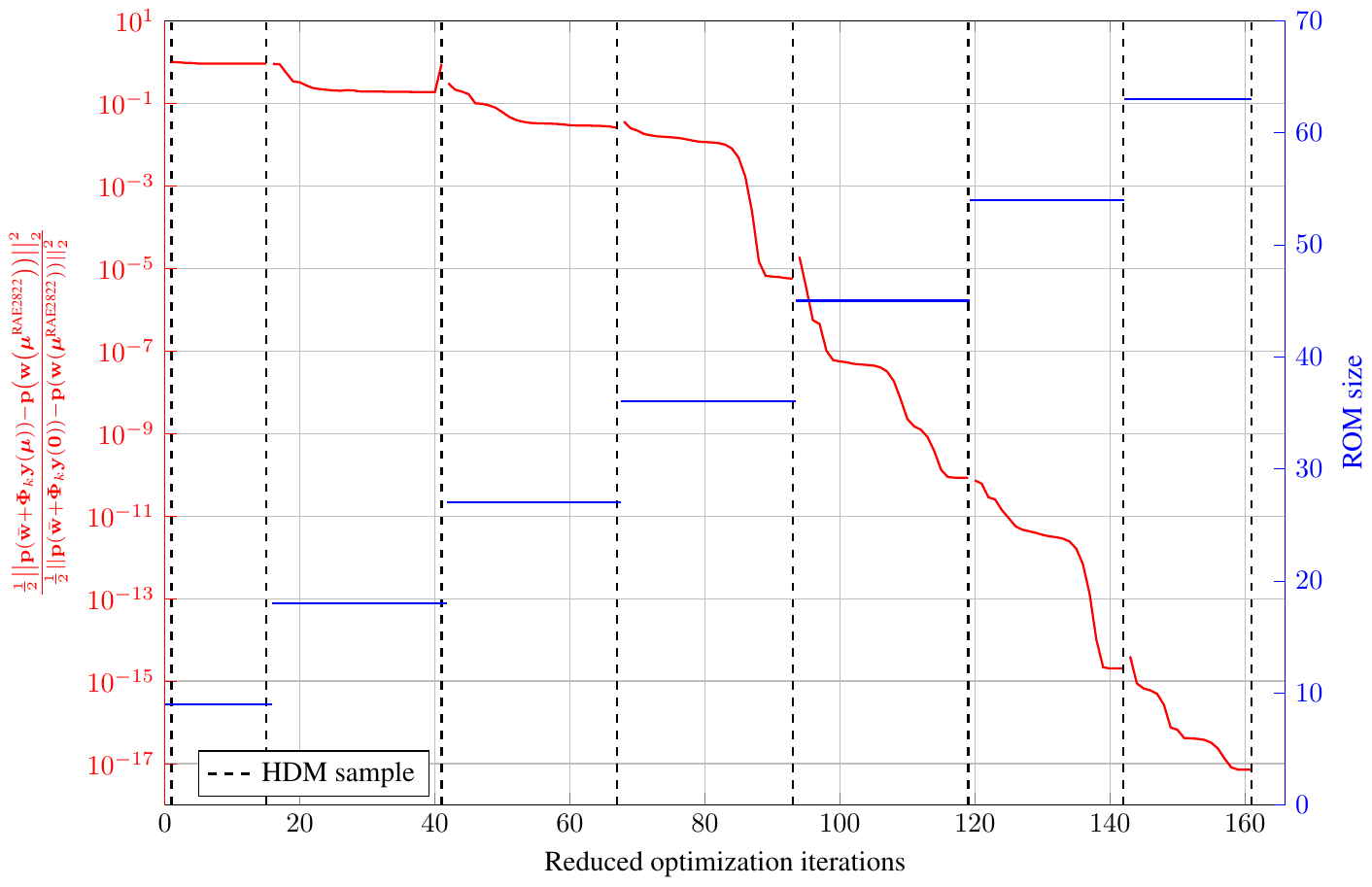}
     \caption{Progression of \emph{reduced} objective function: dashed line indicates an HDM sample and a subsequent update of the ROB}
     \label{fig:match-rom-obj}
\end{figure}

Figure~\ref{fig:match-hdmrom-obj} reports on the evolution of the objective function with the number of performed optimization
iterations, and marks each new HDM query along the optimization trajectory. The reader can observe that the proposed ROM-based 
optimization method performs a total of 160 iterations (Figure~\ref{fig:match-rom-obj}) requiring 346 ROM evaluations (see Table~\ref{tab:summary})
and 7 HDM queries.  From a computational complexity viewpoint, this compares favorably with the 24 HDM-based optimization iterations requiring 29 HDM queries (see Table~\ref{tab:summary}).

Figure~\ref{fig:match-rom-obj} graphically depicts the progression of the \emph{reduced} objective function
across all reduced optimization problems using a dashed line to indicate a new HDM sample and a subsequent update of the ROB.  
For each optimization problem, it also reports the size of the ROB/ROM. 

Figure~\ref{fig:match-rom-fluxnorm} shows the evolution of the HDM residual evaluated at the solution of the ROM --- which
is an indicator of the ROM error --- across all reduced optimization problems, along with the upper bound $\epsilon$, 
defined in (\ref{opt:match-red}) and adapted between reduced optimization problems according to (\ref{eqn:eps-adapt}).  
It is noted here that it is common practice in nonlinear programming software to allow violation of nonlinear constraints during 
an optimization procedure, which explains the residual bound violation seen in this figure.
Figure~\ref{fig:match-rom-fluxnorm} also shows that the ROM solution coincides with the HDM solution at the initial 
condition of each optimization problem, as expected from the choice of the reference vector $\bar\wbold$ discussed 
in Section~\ref{sec:init-guess-rom}.  In the first few reduced optimization problems where the iterate solutions are far 
from the optimal solution, the residual grows as the iterates move into areas of the parameter space away from HDM 
samples. However, near the optimal solution, the residual remains small as the optimization iterates remain in a 
small neighborhood of the most recent HDM sample.

\begin{Rem}
There are two mechanisms that prevent the reduced optimization problem from venturing into regions of the parameter space where it lacks accuracy: \begin{inparaenum}[(1)] \item the objective function, and \item the nonlinear trust region. \end{inparaenum}  
In the present inverse design example, the objective function is mostly sufficient to keep the ROM in regions of accuracy, as can 
be seen from Figure~\ref{fig:match-rom-fluxnorm} where the trust region bound is only reached once and the upper bound always increases.  For other objective functions such as drag or downforce, the nonlinear trust region will be necessary as it is likely that an 
inaccurate ROM can predict a lower value in such objective functions than is actually present.  In practice, inaccurate ROMs have 
been observed to predict the unphysical situation of negative drag (i.e. thrust), which motivates the need for the residual-based trust region. 
\end{Rem}

\begin{savenotes}
\begin{table}[h]
\centering
\begin{tabular}{  >{\centering\arraybackslash}m{1.5in}  >{\centering\arraybackslash}m{1.25in}  >{\centering\arraybackslash}m{1.25in}} \toprule
                       & HDM-based optimization &  ROM-based optimization \\\midrule
\# of HDM Evaluations  &   29             &          7\footnote{The last HDM sample in Figure~\ref{fig:match-hdmrom-obj} was not included in this count as the residual-based error indicator is small at this configuration (Figure~\ref{fig:match-rom-fluxnorm}).  A similar argument could also be made for the $7$th HDM sample.}       \\ 
\# of ROM Evaluations  &   -              &         346       \\ 
\rule[-2.5ex]{0pt}{6ex} $\displaystyle{\frac{\norm{\mubold^* - \mubold^{RAE2822}}}{\norm{\mubold^{RAE2822}}}} $ & $2.28 \times 10^{-3}\%$ & $4.17 \times 10^{-6}\%$\\\bottomrule \addlinespace[3\belowrulesep]
\end{tabular}
\caption{Performance of the HDM- and ROM-based optimization methods} \label{tab:summary}
\end{table}
\end{savenotes}


\section{CONCLUSIONS} \label{sec:conc}

The proposed adaptive approach for using Reduced-Order Models (ROMs) as surrogate models in PDE-constrained optimization 
breaks the traditional offline-online framework of model order reduction.  It builds a Reduced-Order Basis (ROB) using a 
variant of the POD method that incorporates both state and sensitivity snapshots. This is because the low-dimensionality 
assumption of model reduction implicitly assumes that the ROB can approximately represent the state vector and its sensitivities.  
The resulting ROM is used as a surrogate for the HDM in the optimization problem and a nonlinear trust region constraint is 
included to restrict the growth of a residual-based error indicator.  The adaptive approach samples the HDM at the solution of 
this reduced optimization problem and uses the new information to update the ROB and define a new reduced optimization problem.

The minimum-error reduced sensitivity framework introduced in this paper offers the advantage of returning the best approximation 
to the HDM sensitivities that is available in the given ROB, for some norm. In the context of the Least-Squares Petrov-Galerkin 
projection, the minimum-error reduced sensitivities approach the true reduced sensitivities as the error in the ROM approaches zero.  The proposed minimum-error reduced sensitivity framework ensures that the approximation of the sensitivities cannot be degraded 
by state vector snapshots in the ROB.

Most importantly, the proposed progressive approach to ROM-constrained optimization avoids the difficult task of building a single 
{\it global} ROB that is accurate throughout the parameter space. Because it also avoids sampling in regions of the parameter space 
that have not been visited during the optimization procedure, it requires fewer queries to the HDM than traditional alternatives.  
A factor of 4 fewer HDM queries were observed on a problem from aerodynamic shape optimization, where the optimal solution 
was recovered to machine precision. Transforming this computational complexity advantage into a CPU advantage requires however 
equipping this proposed approach for ROM-based optimization with hyperreduction, a step which is the subject of on-going work.

\appendix
\section{Low-Rank SVD Update Algorithms}

\begin{algorithm}
  \caption{Brand's Algorithm for low-rank SVD updates: appending vector} \label{alg:svd-append}
  \begin{algorithmic}[1]
    \REQUIRE Data matrix, $\Xbold \in \Rbb^{m \times n}$ of rank $r$; thin SVD of data matrix, $\Xbold = \Ubold \Sigmabold \Vbold^T$; and full-rank matrices of vectors to append data matrix, $\Ybold \in \Rbb^{m \times k}$. 
    \ENSURE SVD of updated data matrix: $\begin{bmatrix}\Xbold & \Ybold\end{bmatrix} = \bar\Ubold \bar\Sigmabold\bar\Vbold^T$
    \STATE Compute $\Mbold = \Ubold^T\Ybold \in \Rbb^{r \times k}$
    \STATE Compute $\bar\Pbold = \Ybold - \Ubold\Mbold \in \Rbb^{m\times k}$
    \STATE Compute QR decomposition of $\bar\Pbold = \Pbold\Rbold_A$, where $\Pbold \in \Rbb^{m \times k}$, $\Rbold_A \in \Rbb^{k \times k}$
    \STATE Form $\Kbold = \begin{bmatrix} \Sigmabold & \Mbold \\ \zerobold & \Rbold_A\end{bmatrix} \in \Rbb^{(r + k) \times (r + k)}$
    \STATE Compute SVD of $\Kbold = \Cbold \Sbold \Dbold^T$, where $\Cbold, \Sbold, \Dbold \in \Rbb^{(r + k) \times (r + k)}$
    \STATE $\bar\Ubold = \begin{bmatrix} \Ubold & \Pbold\end{bmatrix} \Cbold$
    \STATE $\bar\Sigmabold = \Sbold$
    \STATE $\bar\Vbold = \begin{bmatrix} \Vbold & \zerobold \\ \zerobold & \Ibold\end{bmatrix}\Dbold$
  \end{algorithmic}
\end{algorithm}

\begin{algorithm}
  \caption{Brand's Algorithm for low-rank SVD updates: translating columns} \label{alg:svd-translate}
  \begin{algorithmic}[1]
    \REQUIRE Data matrix, $\Xbold \in \Rbb^{m \times n}$ of rank $r$; thin SVD of data matrix, $\Xbold = \Ubold \Sigmabold \Vbold^T$; and desired translation vector, $\abold \in \Rbb^m$
    \ENSURE SVD of updated data matrix: $\Xbold + \abold\onebold^T = \bar\Ubold \bar\Sigmabold\bar\Vbold^T$
    \STATE Compute $\nbold = \Vbold^T\onebold$, $\qbold = \onebold - \Vbold\nbold$, $q = \norm{\qbold}_2$, and $\Qbold = \frac{1}{q}\qbold$
    \STATE Compute $\mbold = \Ubold^T\abold \in \Rbb^r$
    \STATE Compute $\pbold = \abold - \Ubold\mbold \in \Rbb^m$
    \STATE Define $\hat r \in \Rbb$ and $\vbold \in \Rbb^N$ such that: $\pbold = \hat r\vbold$ where $\norm{\vbold}_2 = 1$
    \STATE Form $\Kbold = \begin{bmatrix} \Sigmabold +\mbold\nbold^T & q\mbold \\ \hat r\nbold & \hat rq\end{bmatrix} \in \Rbb^{(r + 1) \times (r + 1)}$
    \STATE Compute SVD of $\Kbold = \Cbold \Sbold \Dbold^T$, where $\Cbold, \Sbold, \Dbold \in \Rbb^{(r+1) \times (r+1)}$
    \STATE $\bar\Ubold = \begin{bmatrix} \Ubold & \vbold\end{bmatrix} \Cbold$
    \STATE $\bar\Sigmabold = \Sbold$
    \STATE $\bar\Vbold = \begin{bmatrix} \Vbold & \Qbold\end{bmatrix}\Dbold$
  \end{algorithmic}
\end{algorithm}

\bibliographystyle{ieeetr}
\bibliography{progressbasisopt_final}

\ack
The first author would like to thank Kyle Washabaugh for enlightening discussions on aerodynamic optimization 
and the incorporation of sensitivity snapshots in a ROB. He also thanks Kurt Maute for providing access to the SDESIGN software.
The first author also acknowledges the support of the Department of Energy Computational Science Graduate Fellowship. 
The second author acknowledges partial support by the Army Research Laboratory through the Army High Performance Computing 
Research Center under Cooperative Agreement W911NF-07-2-0027, and partial support by the Office of Naval Research under 
Grant N00014-11-1-0707. The content of this publication does not necessarily reflect the position or policy of any of these supporters, and no official endorsement should be inferred.

\end{document}

%% file: mycommands
\newcommand{\beq}{\begin{equation}}
\newcommand{\eeq}{\end{equation}}
\newcommand{\esplit}{\end{split}}
\newcommand{\beqalign}{\begin{array}{rl}}
\newcommand{\eeqalign}{\end{array}}

\theoremstyle{remark}
\newtheorem*{Rem}{Remark} 

\newcommand{\func}[3]{\ensuremath{#1 : #2 \rightarrow #3}}

\newcommand{\norm}[1]{\ensuremath{|| #1 ||}}

\newcommand{\pder}[2]{\ensuremath{\frac{\partial #1}{\partial #2}}} 

\newcommand{\oder}[2]{\ensuremath{\frac{\mathrm{d} #1}{\mathrm{d} #2}}} 

\newcommand{\Dcal}{\ensuremath{\mathcal{D}}}

\newcommand{\Fcal}{\ensuremath{\mathcal{F}}}

\newcommand{\Jcal}{\ensuremath{\mathcal{J}}}

\newcommand{\Lcal}{\ensuremath{\mathcal{L}}}

\newcommand{\Scal}{\ensuremath{\mathcal{S}}}

\newcommand{\Wcal}{\ensuremath{\mathcal{W}}}

\newcommand{\Rbb}{\ensuremath{\mathbb{R} }}

\newcommand\Abold{\ensuremath{\mathbf{A}}}
\newcommand\Bbold{\ensuremath{\mathbf{B}}}
\newcommand\Cbold{\ensuremath{\mathbf{C}}}
\newcommand\Dbold{\ensuremath{\mathbf{D}}}

\newcommand\Ibold{\ensuremath{\mathbf{I}}}

\newcommand\Kbold{\ensuremath{\mathbf{K}}}

\newcommand\Mbold{\ensuremath{\mathbf{M}}}

\newcommand\Pbold{\ensuremath{\mathbf{P}}}
\newcommand\Qbold{\ensuremath{\mathbf{Q}}}
\newcommand\Rbold{\ensuremath{\mathbf{R}}}
\newcommand\Sbold{\ensuremath{\mathbf{S}}}

\newcommand\Ubold{\ensuremath{\mathbf{U}}}
\newcommand\Vbold{\ensuremath{\mathbf{V}}}

\newcommand\Xbold{\ensuremath{\mathbf{X}}}
\newcommand\Ybold{\ensuremath{\mathbf{Y}}}

\newcommand\abold{\ensuremath{\mathbf{a}}}
\newcommand\bbold{\ensuremath{\mathbf{b}}}
\newcommand\cbold{\ensuremath{\mathbf{c}}}

\newcommand\ebold{\ensuremath{\mathbf{e}}}

\newcommand\mbold{\ensuremath{\mathbf{m}}}
\newcommand\nbold{\ensuremath{\mathbf{n}}}
\newcommand\pbold{\ensuremath{\mathbf{p}}}
\newcommand\qbold{\ensuremath{\mathbf{q}}}

\newcommand\ubold{\ensuremath{\mathbf{u}}}
\newcommand\vbold{\ensuremath{\mathbf{v}}}
\newcommand\wbold{\ensuremath{\mathbf{w}}}
\newcommand\xbold{\ensuremath{\mathbf{x}}}
\newcommand\ybold{\ensuremath{\mathbf{y}}}

\newcommand\mubold{\ensuremath{\boldsymbol{\mu}}}

\newcommand\nubold{\ensuremath{\boldsymbol{\nu}}}

\newcommand\Phibold{\ensuremath{\boldsymbol{\Phi}}}
\newcommand\Sigmabold{\ensuremath{\boldsymbol{\Sigma}}}

\newcommand\Thetabold{\ensuremath{\boldsymbol{\Theta}}}

\newcommand\Psibold{\ensuremath{\boldsymbol{\Psi}}}

\newcommand\zerobold{\ensuremath{\mathbf{0}}}
\newcommand\onebold{\ensuremath{\mathbf{1}}}